\documentclass[11pt]{amsart}
\usepackage[leqno]{amsmath}   
\usepackage{amsthm}   
\usepackage{amssymb}  
\usepackage{amscd}   
\usepackage[matrix,arrow]{xy}
\usepackage{verbatim}
\def\wh#1{\widehat{#1}}
\def\wt#1{\widetilde{#1}}
\theoremstyle{plain}
    \newtheorem{theorem}{Theorem}[subsection]
    
    \newtheorem{proposition}[theorem]{Proposition}
    \newtheorem{lemma}[theorem]{Lemma}
    \newtheorem{corollary}[theorem]{Corollary}
\theoremstyle{definition}
    \newtheorem{definition}[theorem]{Definition}
    
    \newtheorem{example}[theorem]{Example}
    \newtheorem{remark}[theorem]{Remark}
\def\Alphabet{A,B,C,D,E,F,G,H,I,J,K,L,M,N,O,P,Q,R,S,T,U,V,W,X,Y,Z}
\def\alphabet{a,b,c,d,e,f,g,h,i,j,k,l,m,n,o,p,q,r,s,t,u,v,w,x,y,z}
\def\endpiece{xxx}
\def\makeAlphabet[#1]{\expandafter\makeA#1,xxx,}
\def\makealphabet[#1]{\expandafter\makea#1,xxx,}
\def\makeA#1,{\def\temp{#1}\ifx\temp\endpiece\else%
\mkbb{#1}\mkfrak{#1}%
\mkcal{#1}\expandafter\makeA\fi}%
\def\makea#1,{\def\temp{#1}\ifx\temp\endpiece\else\mkfrak{#1}
\expandafter\makea\fi}%
\def\mkbb#1{\expandafter\def\csname bb#1\endcsname{\mathbb{#1}}}
\def\mkfrak#1{\expandafter\def\csname fr#1\endcsname{\mathfrak{#1}}}
\def\mkcal#1{\expandafter\def\csname c#1\endcsname{\mathcal{#1}}}
\def\makeop[#1]{\xmakeop#1,xxx,}
\def\mkop#1{\expandafter\def\csname #1\endcsname{{\mathrm{#1}}}} %
\def\xmakeop#1,{\def\temp{#1}\ifx\temp\endpiece\else\mkop{#1}\expandafter\xmakeop\fi}%
\makeAlphabet[\Alphabet]
\makealphabet[\alphabet]
\makeop[Alt,End,Ext,Hom,Sym,Tor]
\makeop[CH,Pic,div]
\makeop[Ker,Im,Coker,Coim,id,pr,sgn]
\makeop[Spec,Spf,Spm]
\makeop[Re,Im]%
\makeop[dR,Nis,crys,rig,syn,tor,Cont]
\makeop[Gal,Res,ab]
\makeop[pol]
\makeop[Var,an,Isoc,univ,as]%
\makeop[Gl,Sl]
\makeop[Lie,Ala,Al]
\makeop[tay,tr,inf]
\makeop[Sat,Der]

\def\BGinf{B_\cdot^{1}G}
\def\EGinf{E_\cdot^{1}G}

\def\isom{\cong}
\def\verk{\circ}

\def\Gm{\bbG_m}

\renewcommand{\Bbb}{\mathbb}  

\newcommand{\Q}{{\Bbb{Q}}}  
\newcommand{\C}{{\Bbb{C}}}  
\newcommand{\Z}{{\Bbb{Z}}}  
\newcommand{\A}{{\Bbb{A}}}  

\newcommand{\F}{{\Bbb{F}}}  

\newcommand{\Xf}{\mathfrak{X}}
\newcommand{\Yf}{\mathfrak{Y}}
\newcommand{\Pf}{\mathfrak{P}}

\newcommand{\Oh}{\mathcal{O}}

\newcommand{\tensor}{\otimes}   

\newcommand{\indlim}{\varinjlim}

\newcommand{\bew}{\begin{proof}}
\newcommand{\bewende}{\end{proof}}
\theoremstyle{plain}
\newtheorem{prop}[theorem]{Proposition}
\newtheorem{thm}[theorem]{Theorem}

\newtheorem{cor}[theorem]{Corollary}

\theoremstyle{definition}
\newtheorem{defn}[theorem]{Definition}
\newtheorem{rem}[theorem]{Remark}
\newtheorem{ex}[theorem]{Example}

\newcommand{\g}{\mathfrak{g}}
\newcommand{\gl}{\mathfrak{gl}}
\newcommand{\GL}{\mathrm{GL}}

\newcommand{\cont}{\mathrm{cont}}
\newcommand{\BK}{\mathrm{BK}}
\newcommand{\DR}{\mathrm{DR}}

\newcommand{\et}{\mathrm{et}}
\newcommand{\cone}{\mathrm{Cone}}
\newcommand{\Fil}{\mathrm{Fil}}
\newcommand{\la}{\mathrm{la}}
\renewcommand{\sp}{\mathrm{sp}}
\newcommand{\maps}{\mathrm{Maps}}
\newcommand{\Sp}{\mathrm{Sp}}
\newcommand{\borel}{b_p}
\newcommand{\BG}{B.G}
\newcommand{\BGLN}{B.\GL_N}
\newcommand{\BcG}{B.\cG}
\newcommand{\alg}{\mathrm{alg}}
\newcommand{\eval}{\mathrm{ev}}

\newcommand{\Ohla}[1]{\Oh({#1}^\la)}
\newcommand{\Ohdagger}[1]{\Oh({#1}^\dagger)}
\newcommand{\Ohcont}[1]{\Oh({#1}^\cont)}
\newcommand{\Ohdelta}[1]{\Oh({#1}^\delta)}
\newcommand{\Ohh}[1]{\Oh(#1)}
\newcommand{\Ohalg}[1]{\Oh({#1}^\alg)}

\newcommand{\Ohlatwo}[2]{\Oh({#1}^\la,#2)}

\newcommand{\Ohconttwo}[2]{\Oh({#1}^\cont,#2)}
\newcommand{\Ohdeltatwo}[2]{\Oh({#1}^\delta,#2)}

\begin{document}
\title[A $p$-adic Borel regulator and Bloch-Kato exponential]{A $p$-adic analogue of the Borel regulator and the Bloch-Kato
exponential map}
\author{Annette Huber}
\address{Annette Huber, Math. Institut, Universit\"at Leipzig, Johannisgasse 26, 04109 Leipzig, Germany}
\email{huber@mathematik.uni-leipzig.de}
\author{Guido Kings}
\address{Guido Kings, NWF I - Mathematik,   Universit\"at Regensburg,  93040 Regensburg, Germany}
\email{guido.kings@mathematik.uni-regensburg.de}
\begin{abstract}In this paper we define a $p$-adic analogue of the Borel regulator for the
$K$-theory of $p$-adic fields. The van Est isomorphism in the construction of the
classical Borel regulator is replaced by the Lazard isomorphism. The main
result relates this $p$-adic regulator to the Bloch-Kato exponential and the
Soul\'e regulator. On the way we give a new description of the Lazard isomorphism
for certain formal groups. 
\end{abstract}
\date{\today} 
\maketitle
\tableofcontents

\section*{Introduction}
The classical Borel regulator 
$$
	b_\infty:K_{2n-1}(\bbC)\to \bbC
$$
plays a decisive role in the study  of
algebraic number fields. Its simple description makes it possible to relate 
this regulator quite directly to special values of zeta functions of number fields,
an insight we owe to Borel \cite{Bo}. Later, when Beilinson formulated his famous 
conjectures, he was able to relate his regulator $r_\infty$ to Borel's regulator map,
so that the computations by Borel implied one of the most impressive confirmations
of Beilinson's conjectures.

The aim of this paper is to  study an analogue of the Borel regulator for local fields
and to relate it to the Soul\'e regulator via the exponential map of Bloch and Kato.

Let $K/\bbQ_p$ be a finite extension and $R\subset K$ its valuation ring.
Soul\'e has defined a regulator map
$$
	r_p:K_{2n-1}(K)\to H^1_\et(K,\bbQ_p(n))
$$
between the $K$-theory of $K$ and \'etale cohomology.
We define in Section \ref{padicborel}, in complete analogy with
the classical Borel regulator for $\bbC$, a map
$$
	b_p:K_{2n-1}(R)\to K
$$
such that the diagram
\begin{equation}\label{introdiag}
\begin{xy}
\xymatrix{K_{2n-1}(R)\ar[rr]^{r_p}\ar[rrd]^{b_p}&&H^1_\et(K,\bbQ_p(n))\\
&&K=D_\dR(\bbQ_p(n))\ar[u]_{\exp_{\BK}}}
\end{xy}
\end{equation}
commutes for all $n\ge 1$. Here 
$ \exp_{\BK}:K\to H^1_\et(K,\bbQ_p(n))$ is the Bloch-Kato exponential map, which is an isomorphism for $n>1$.

The idea to define $b_p$ follows closely Borel's construction (as reviewed in Section \ref{classicborel}), 
only that we replace the van Est isomorphism
$$
	H^i(\frg,\frk,\bbC)\isom H^i_\cont(G(\bbC),\bbC)
$$
between relative Lie algebra cohomology and continuous group cohomology, by the
Lazard isomorphism 
$$
	H^i(\frg,\bbQ_p)\isom H^i_\cont(G(R),\bbQ_p)
$$
(for details see Section \ref{padicborel}). \\

The commutative diagram (\ref{introdiag}) above can be seen as a natural generalisation of the explicit
reciprocity law for $\Gm$
$$
\begin{xy}
\xymatrix{R^*\ar[rrd]^{\log_p}\ar[rr]^{\partial}&&H^1_\et(K,\bbQ_p(1))\\
&&K\ar[u]_{\exp_{\BK}},}
\end{xy}
$$
(where $\partial$ is the Kummer map)
established by Bloch and Kato in their work on the exponential map \cite{BK} 3.10.1.\\

Karoubi has defined a $p$-adic regulator for $p$-adic Banach algebras in \cite{Karoubi},
which was studied in more detail in the thesis of Hamida \cite{Hamida}. It would
be very interesting to investigate the exact relation of her results with our construction.\\

On the way to establish the diagram (\ref{introdiag}), we also get a new description
of the Lazard isomorphism. Namely,
we show that for a smooth algebraic group scheme $H/R$ with formal group $\cH$ the
Lazard isomorphism $H^i_\la(\cH,\bbQ_p)\isom H^i(\Lie H,\bbQ_p)$ is induced by the map
\begin{align*}
	\Phi:\Ohla{\cH}^{\otimes n}&\to {\bigwedge}^n\Lie H\\
	f_1\otimes\ldots\otimes f_n& \mapsto df_1\wedge\ldots\wedge df_n
\end{align*}
familiar from cyclic homology.\\

Unfortunately, the proof of diagram (\ref{introdiag}) is quite involved and 
technical. The idea is to follow Beilinson's proof, which leads to a comparison
between the classical Borel regulator and Beilinson's regulator in Deligne cohomology.
Our strategy adapts Beilinson's ideas to the $p$-adic case replacing Deligne cohomology
by syntomic cohomology.\\

The paper is organised as follows: in the first section we construct 
the $p$-adic analogue of the Borel regulator and formulate our main results.
We also draw some immediate consequences.

The second section recalls what we need about rigid syntomic cohomology and introduces
the evaluation map \ref{eval}, which is the main tool in the comparison with the locally analytic
group cohomology. 

The third section relates this evaluation map to locally analytic group cohomology via
the suspension map and Lie algebra cohomology. Here we follow quite closely Beilinson's
ideas as explained by Rapoport \cite{Ra} and Burgos \cite{Burgos}.

The fourth section shows finally that the map between locally analytic group cohomology
and Lie algebra cohomology is indeed the Lazard isomorphism.

The final section collects the loose ends and finishes the proof.\\

It is a great pleasure to thank Jos\'e Burgos for answering some questions about
Beilinson's proof of his comparison results and Elmar Gro{\ss}e-Kl\"onne for his help
with rigid cohomology. Of course the whole article owes its very existence to the 
beautiful ideas of Borel and Beilinson.

\section*{Notation}
In this section we collect various notations, which will be needed later. 
\subsection{Rings of functions}
Let $p$ be a fixed prime. Let $R$ be a discrete valuation ring finite over $\Z_p$ with uniformiser $\pi$, residue field $k$ and field of fractions $K$.  
Throughout let $G$ be $\GL_N$,  the general linear group over $R$ and $\g$ its $K$-Lie algebra, ie. $\gl_N=M_N(K)$ the
$N\times N$-matrices.

To  a smooth $R$-scheme $X$, we attach a number of spaces. The set $X(R)$ is denoted $X^\delta$. It
carries a natural structure of (locally) analytic manifold over $K$, which is denoted $X^\la$.
The underlying topological space of $X^\la$ is denoted $X^\cont$.
In Section \ref{syn} we are also going to consider its structure as overconvergent rigid analytic manifold, denoted $X^\dagger$ (see Example \ref{schemedagger}). In the case $X=G$, the set
$G^\delta$ is a group, $G^\cont$ a topological group and $G^\la$ is a $K$-Lie group.

We also denote
\begin{itemize}
\item $\Ohh{X}$ the global sections of the $R$-scheme $X$,
\item $\Ohalg{X}$ the global sections of the $K$-scheme $X\times_R K$,
\item $\Ohdelta{X}$ the ring of set-theoretic $K$-valued functions on $X(R)$,
\item $\Ohcont{X}$ the ring of continuous $K$-valued functions on $X(R)$,
\item $\Ohla{X}$ the ring of $K$-valued locally analytic functions on $X^\la$, i.e. functions which are locally representable by convergent power series with coefficients in $K$,
\item $\Ohdagger{X}$ the ring of overconvergent rigid analytic functions on 
$X^\dagger$.
\end{itemize}
If $A$ is a $K$-vector space, we denote
$\Ohlatwo{X}{A}$, $ \Ohdeltatwo{X}{A}$ etc. the corresponding functions with values in $K$.

If $X.$ is a smooth simplicial $R$-scheme, $\Ohh{X.}$ etc. are cosimplicial rings.

\subsection{The classifying space}

We collect some standard material on the classifying space $B.H$ thereby fixing our notations.

Let $H/R$ be a reductive algebraic group over a ring $R$ or a formal group.
As usual, let $E.H$ be the simplicial space with
$$
	E_nH=H\times\dots\times H\mbox{ $n+1$-times}
$$
with the usual face and degeneracy maps. The group $H$ acts on $E.H$ on
the right via
$$
	(h_0,\dots, h_n)h:=(h_0h,\dots, h_nh).
$$
The quotient by this action is the classifying space $B.H$ and has the explicit description
$$
	B_nH=H\times\dots\times H\mbox{ $n$-times}
$$
with degeneracies $\sigma^i(h_1,\dots, h_n)=(h_1,\dots, h_i,1,h_{i+1},\dots, h_n)$ and faces
\begin{align*}
	\delta^0(h_1,\dots, h_n)&=(h_2,\dots, h_n)\\
	\delta^i(h_1,\dots, h_n)&=(h_1,\dots, h_ih_{i+1},\dots ,h_n)\mbox{ for }i=1,\dots, n-1\\
	\delta^n(h_1,\dots, h_n)&=(h_1,\dots, h_{n-1}).
\end{align*}
The map $E.H\to B.H$ sends $(h_0,\dots, h_n)\mapsto (h_0h_1^{-1},\dots, h_{n-1}h_n^{-1})$.

\subsection{Group cohomology}
Let $H$ be an abstract group.
Let $F^*$ be a complex of $\Q_p$-vector spaces   with trivial operation of $H$.
Let
$H^i(H,F^*)$ be group hypercohomology of $H$ with coefficients in $F$, i.e.,
in the notation introduced above, cohomology of the complex associated to the cosimplicial complex
$\Ohdeltatwo{B.H}{F^*}$. 

For a $K$-Lie group $H$  and a complex $F^*$ of $K$-vector spaces, 
let $H^i_\cont(H,F^*)$ be {\em continuous group cohomology}, i.e., cohomology of 
$\Ohconttwo{B.H}{F^*}$ and 
$H^i_\la(H,F^*)$ the {\em locally analytic group cohomology}, i.e., the  cohomology of 
$\Ohlatwo{B.H}{F^*}$.
The inclusions of function spaces induce natural transformations
\[ H^i_\la(H,F^*)\to H^i_\cont(H,F^*)\to H^i(H,F^*).\]
If $H$ is a smooth algebraic group over $R$, we abbreviate
$H^i_\cont(H,F^*)=H^i_\cont(H^\cont,F^*)$ etc.

\subsection{Lie algebras and their cohomology}\label{lie}
\begin{defn}\label{lie-algebra}Let $H$ be a smooth algebraic group over $R$.
We denote
\[ \Lie H=\Der_K(\Ohalg{H}_e,K)\]
its algebraic tangent space at $e$. This is the $K$-{\em Lie algebra} of $H$.
\end{defn}
\begin{example} For $G=\GL_N$, we have naturally $\Lie G\isom \g$.
\end{example}

By definition, there is a natural pairing
\[ \Ohalg{H}_e\times \Lie H\to K\]
It extends to locally analytic functions and induces
\begin{align*} \Ohla{H}&\to\Lie H^\vee\\
                  f&\mapsto df(e)
\end{align*}

\begin{defn}Let $\frh$ be a Lie algebra over $K$. Then {\em Lie algebra cohomology}
$H^i(\frh,K)$ is defined as cohomology of the complex 
$\bigwedge^*\frh^\vee_K$ with differential induced by the dual of the 
Lie bracket $\frh\tensor\frh\to \frh$. 
\end{defn}

\begin{ex}For $\g=\gl_N$, $H^*(\g,K)$ is an exterior algebra on primitive elements $p_n\in H^{2n-1}(\g,K)$ for $n=1,\dots, N$.
\end{ex}
We need a precise normalisation of these element and follow 
\cite{Burgos} Example 5.37.

\begin{defn}\label{primitive} Let $\g=\gl_N$. For $n\leq N$ let $p_n\in H^{2n-1}(\g,K)$ be the map
\[ p_n:{\bigwedge}^{2n-1}\g\to K\]
which attaches to $(x_1,\dots, x_{2n-1})\in\g^{\oplus 2n-1}$ the value
\[ \frac{((n-1)!)^2}{(2n-1)!}\sum_{\sigma\in \mathfrak{S}_{2n-1}}\operatorname{sgn}(\sigma)\operatorname{Tr}(x_1\circ\dots\circ x_{2n-1})\]
where $\mathfrak{S}_{2n-1}$ is the symmetric group, $\operatorname{Tr}$ is the
trace map from $\gl_N$ to $K$ and $\circ$ is matrix multiplication.
\end{defn}
\begin{rem}
	The exact form of this primitive element $p_n$ is not necessary for our main result. What we need is  the
	element $p_n\in H^{2n-1}(\g,K)$ which is the image of the Chern class $c_n\in H_{\dR}^{2n}(B.\GL_N/K)$ under 
	the suspension map $s_G:H_{\dR}^{2n}(B.\GL_N/K)\to H_{\dR}^{2n-1}(\GL_N/K)\isom H^{2n-1}(\g,K)$ (see 
	Section~\ref{suspdefn} for more on the suspension map).
\end{rem}

\section{Construction of the $p$-adic Borel regulator and statement of the main results}\label{setup}

\subsection{Review of the classical Borel regulator}\label{classicborel}
We briefly review the construction of the classical Borel regulator
$$
	b_{\infty}:K_{2n-1}(\bbC)\to \bbC
$$
from \cite{Bor}.
Recall that the $K$-group $K_{2n-1}(\bbC)$ maps via the Hurewicz map to the group homology
$H_{2n-1}(\GL(\bbC),\bbC)$ of the infinite linear group $\GL(\bbC):=\indlim_N\GL_N(\bbC)$. 
Thus one needs to define a system of compatible maps for all big enough $N$ (also called $b_\infty$)
$$
	b_\infty:H_{2n-1}(\GL_N(\bbC),\bbC)\to \bbC,
$$
or, by duality for group cohomology, elements
$$
	b_\infty\in H^{2n-1}(\GL_N(\bbC),\bbC),
$$
compatible with enlarging $N$. In fact, using stabilisation,
these elements live in $H_{2n-1}(\GL_N(\bbC),\bbC)$ for fixed $n$ and $N$ big enough.
Borel constructs this $b_\infty$ with the \emph{van Est isomorphism} for $\GL_N(\bbC)$ considered as
real Lie group
$$
	H^{2n-1}(\frg\frl_N,\fru_N,\bbC)\isom H^{2n-1}_\cont(\GL_N(\bbC),\bbC)
$$
between relative Lie algebra cohomology for $\frg\frl_N$ with respect to the Lie algebra
$\fru_N$ of the unitary group $U_N$ and continuous group cohomology.
To compute $H^{2n-1}(\frg\frl_N,\fru_N,\bbC)$ observe that the compact form of
$\GL_N(\bbC)$ is $U_N\times U_N$ so that we have isomorphisms
$$
H^{2n-1}(\frg\frl_N,\fru_N,\bbC)\isom H^{2n-1}(\fru_N\oplus \fru_N,\fru_N,\bbC)\isom H^{2n-1}(\fru_N,\bbC)\isom H^{2n-1}(\frg\frl_N,\bbC).
$$
Recall from \ref{primitive} that we have defined a primitive element $p_n\in H^{2n-1}(\frg\frl_N,\bbC)$.
\begin{defn}
	The \emph{Borel regulator} $b_\infty$ is the image of $p_n$ under the composition
	$$
		H^{2n-1}(\frg\frl_N,\bbC)\isom H^{2n-1}(\frg\frl_N,\fru_N,\bbC)\isom H^{2n-1}_\cont(\GL_N(\bbC),\bbC)\to 
		H^{2n-1}(\GL_N(\bbC),\bbC).
	$$
\end{defn}

\begin{rem}
	Beilinson shows in \cite{Be} (see also \cite{Ra}, \cite{Burgos}) that under the 
	identification of Deligne cohomology $H^{1}_{\cD}(\Spec\bbC,\bbR(n))\isom\bbC$
	his regulator $r_\infty$ coincides with Borel's regulator $b_\infty$ up to a rational number. It was later shown by Burgos \cite{Burgos} that in fact $2r_\infty=b_\infty$. 
\end{rem}

\subsection{A $p$-adic analogue of Borel's regulator}\label{padicborel}
Recall that $K/\bbQ_p$ is a finite extension with valuation ring $R$.
Our aim is to define a $p$-adic analogue of Borel's regulator. We will construct
a map
\[ \borel:K_{2n-1}(R)\to K.\]
As before, we have the Hurewicz map
$K_{2n-1}(R)\to  H_{2n-1}(\GL(R),K)$
and 
it suffices to define a compatible system of maps
\[\borel:H_{2n-1}(\GL_N(R),K)\to K\hspace{3ex}\text{for $N\geq n$}\]
or, equivalently, elements 
\[ \borel\in H^{2n-1}(\GL_N(R),K).\]
To define $\borel$ we replace the van Est isomorphism in Borel's construction with
a slight generalisation of the \emph{Lazard isomorphism}: consider $\GL_N(R)$ as
a $K$-analytic group and let $\frg\frl_N$ be its $K$-Lie algebra. We will prove:
\begin{thm}[\cite{L}, see also Section \ref{proofsection}]\label{lazardiso}
For all $i\ge 0$ one has an isomorphism
\[ H^i_\la(\GL_N(R),K)\isom H^i(\frg\frl_N,K)\]
between the locally analytic group cohomology and the Lie algebra cohomology.
\end{thm}
\begin{rem} In the case $R=\Z_p$, this is the combination of two isomorphisms shown by Lazard
$$
	H^i_\la(\GL_N(\bbZ_p),\bbQ_p)\isom H^i_\cont(\GL_N(\bbZ_p),\bbQ_p)
$$
and 
$$
	H^i_\cont(\GL_N(\bbZ_p),\bbQ_p)\isom H^i(\frg\frl_N,\bbQ_p),
$$
(see \cite{L} chapter V). For the argument in the general case see Section \ref{proofsection}.
\end{rem}
Now we can define the $p$-adic analogue of the Borel regulator
using again the primitive element $p_n\in H^i(\frg\frl_N,K)$ from Definition \ref{primitive}.
\begin{definition}
	The \emph{$p$-adic Borel regulator} 
	$$
		\borel:K_{2n-1}(R)\to K
	$$
	for $1\le n$, is defined by the element $\borel\in H^{2n-1}(\GL_N(R),K)$ (for $N$ big enough), 
	which is the image of $p_n$ under the
	composition
	$$
		H^i(\frg\frl_N,K)\isom H^i_\la(\GL_N(R),K)\to H^{2n-1}(\GL_N(R),K).
	$$
\end{definition}
\begin{rem}\begin{enumerate}
\item Soul\'e \cite{So2} was the first to study a $p$-adic regulator for $K_{2n-1}(K)$ with values in $\bbZ_p^{[K:\bbQ_p]}$. 
His regulator is defined via Iwasawa theory.
\item Karoubi (see e.g. \cite{Karoubi}) has also defined a $p$-adic regulator with values in
	topological cyclic homology. This was further studied in the thesis of Hamida \cite{Hamida}.
\end{enumerate}
	
\end{rem}

\subsection{The main result}
Before we can formulate our main result, we need to recall the Soul\'e regulator
and the Bloch-Kato exponential map.

Soul\'e \cite{Soule} has defined \emph{regulators} for $n>0$
\[ r_p:K_{2n-1}(K)\to H^1_\et(K,\Q_p(n)),\]
which are just the Chern classes
\[ c_n\in H^{2n-1}(\GL_N(K),H^1_\et(K,\Q_p(n)))\]
induced by the universal Chern classes (via restriction to $B_\cdot GL_N^\delta$)
\[ c_n\in H^{2n}_\et(B_\cdot GL_N,\Q_p(n)).\]
By abuse of notation, we also let $r_p$ be the composition
$$
	 r_p:K_{2n-1}(R)\to K_{2n-1}(K)\to H^1_\et(K,\Q_p(n)).
$$
This $r_p$ is given by the restriction of the $c_n$ above to $H^{2n-1}(\GL_N(R),H^1_\et(K,\Q_p(n)))$.
\begin{rem}
	Note that for $n>1$ the morphism of $\bbQ_p$-vector spaces
	$$
		K_{2n-1}(R)\otimes\bbQ_p \to K_{2n-1}(K)\otimes\bbQ_p
	$$
	is an isomorphism. This follows from the localisation sequence for
	$K$-theory and from Quillen's result that $K_i$ of a finite field is
	torsion for $i\ge 1$.
\end{rem}
Recall (\cite{BK} Definition 3.10) that the \emph{Bloch-Kato exponential} 
\[ \exp_\BK:K\to H^1_\et(K,\Q_p(n))\]
is the connecting morphism of the short exact sequence of continuous $\Gal(\bar{K}/K)$-modules (\cite{BK} Proposition 1.17)
\[ 0\to\Q_p\to B_\crys^{f=1}\oplus B_\DR^+\to B_\DR\to 0\]
tensored with $\Q_p(n)$.
It is an isomorphism for $n>1$. 

Our main result is the following:
\begin{thm}\label{mainthm} 
The following diagram commutes for all $n>0$:
\begin{equation}
\begin{xy}
\xymatrix{K_{2n-1}(R)\ar[rr]^{r_p}\ar[rrd]^{\borel}&&H^1_\et(K,\bbQ_p(n))\\
&&K\ar[u]_{\exp_{\BK}}}
\end{xy}
\end{equation}
\end{thm}
\begin{proof}[Idea of proof:]
\begin{enumerate}
\item By construction, the statement of the theorem can be formulated as follows:
the $p$-adic Borel regulator
$\borel\in H^{2n-1}(\GL_N(R),K)$ is mapped under the  Bloch-Kato exponential 
$\exp_\BK: K\to H^1_\et(K,\Q_p(n))$ to 
the \'etale Chern class $c_n\in H^{2n-1}(\GL_N(R),H^1_\et(K,\Q_p(n)))
$
\item In Proposition \ref{synet} together with Proposition \ref{etale}
we show that $\exp_{BK}\circ c_n$ is the syntomic Chern class. 
\item Using the evaluation map from syntomic cohomology (Theorem \ref{evalthm}) we consider the image of this
syntomic Chern class in analytic group cohomology.
\item We apply the Lazard morphism to this element and show that it can be identified with the primitive element in Lie algebra
cohomology. For this we have to  relate the Lazard morphism to the suspension map (see Theorem \ref{suspdiagram})
\item This finishes the proof because it is well known that the Chern class (in de Rham cohomology) is
mapped to the primitive element under the suspension map. In fact, this is one of the possible definitions of Chern classes.
\end{enumerate}
The detailed proof will be given in Section \ref{proofsection}.
\end{proof}
\begin{rem} One can see the main theorem as a kind of generalised explicit reciprocity law for 
the formal group of $\Gm$. More precisely, the main theorem generalises the diagram
$$
	\begin{xy}
	\xymatrix{R^*\ar[rrd]_{\log_p}\ar[rr]^{\partial}&&H^1_\et(K,\bbQ_p(1))\\
	&&K\ar[u]_{\exp_{\BK}},}
	\end{xy}
$$
(where $\partial$ is the Kummer map)
established by Bloch and Kato in their work on the exponential map \cite{BK} 3.10.1.
\end{rem}
\begin{cor}
\begin{enumerate}
\item  The Soul\'e regulator $c_n\in H^{2n-1}(\GL_N(R),H^1_\et(K,\Q_p(n)))$ is 
the image of an element 
\[ c_n^\cont\in H^{2n-1}_\cont(\GL_N(R),H^1_\et(K,\Q_p(n)))\]
in continuous group cohomology for $N$ big enough.
\item 
 There is an element
\[ (\exp_\BK^*c_n)^\la \in  H^{2n-1}_\la(\GL_N(R),K)\]
in locally analytic group cohomology for $N$ big enough, such that
the composition with the Bloch-Kato-exponential map $K\to H^1_\et(K,\Q_p(n))$
gives $c_n^\cont$.
\end{enumerate}
\end{cor}
\begin{proof} After identifying $r_p$ (i.e, the \'etale Chern class) with $\borel$ this follows from the 
definition of $\borel$. In fact, however, we are going to prove this directly as a step in the proof of the Main Theorem, see Corollary \ref{chernanalytic}.
\end{proof}
\begin{rem}
	In her thesis \cite{Hamida} Hamida also obtains the result that Karoubi's regulator is induced from 
	a continuous (hence also from a locally analytic) group cohomology class.
\end{rem}
\subsection{The Lazard isomorphism as a Taylor series expansion}

In our proof of the main theorem we will give a new description of the 
Lazard isomorphism at least for formal groups associated to smooth algebraic groups.

Let $H/\bbZ_p$ be a smooth linear algebraic group and $\cH$ a
 $p$-saturated group of finite rank (see \cite{L} III Definition 2.1.3 and 2.1.6) with valuation $\omega$, which is an open subgroup of $H(\Z_p)$. The group $\cH$ will always be considered as a $\bbQ_p$-analytic manifold.
Our main example is $H=\GL_N$, $\cH=1+pM_N(\Z_p)$ with $\omega$ the $\inf$-valuation.

 We show in \ref{liealgident} that the Lie algebra $\cL^*$ of $\cH$
in the sense of Lazard (cf. Definition \ref{LazardLie}) can be identified with the algebraic Lie algebra
$$
	\frh\isom \cL^*\otimes \bbQ_p.
$$
Lazard (\cite{L} chapter V) shows in his paper that one has a chain of isomorphisms
$$
	 H^i_\la(\cH,\bbQ_p)\isom H^i_\cont(\cH,\bbQ_p)\isom H^i(\cL^*,\bbQ_p)
$$
where the last isomorphism is induced by an isomorphism between the saturations of 
the continuous group algebra of $\cH$ and the universal enveloping algebra of $\cL^*$.
We give a much simpler description of this isomorphism:
\begin{definition}\label{Phidefn}
	Define a map of complexes
	\begin{align*}
		\Phi:\Ohla{B_n\cH}\isom \Ohla{\cH}^{\hat{\otimes}n}&\to \bigwedge^n\frh^\vee\\
		f_1\tensor\dots\tensor f_n&\mapsto df_1(e)\wedge\dots\wedge df_n(e).
	\end{align*}
\end{definition}
\begin{thm}[see Proposition \ref{taylorent}] Let us identify $\frh\isom \cL^*\otimes \bbQ_p$ so that
$H^i(\cL^*,\bbQ_p)\isom H^i(\frh,\bbQ_p)$.
The Lazard isomorphism coincides  with the map, which is induced by $\Phi$ on cohomology:
$$
\Phi:H^i_\la(\cH,\bbQ_p)\isom H^i(\frh,\bbQ_p).
$$
\end{thm}

\begin{rem}
	The map used here is well-known in connection with cyclic homology and the 
	Hochschild-Kostant-Rosenberg theorem, see \cite{Lod} 1.3.14.
	It would be interesting to study this relation further.
\end{rem}

\section{Syntomic cohomology}\label{syn}

As before let $R$ be a discrete valuation ring finite over $\Z_p$ with field of fractions $K$ and
residue field $k\isom \F_q$. Let $K_0$ be the maximal absolutely unramified subfield of $K$ and 
$R_0\subset R$ its ring of integers.
There are no conditions on ramification. 

 Let $X$ be a smooth
scheme over $R$. We are going to review the definition of (a certain version of) syntomic cohomology and prove properties of the syntomic Chern classes on $K_i(R)$. The construction follows Besser's $\tilde{H}^i_\syn(X,n)$ \cite{Be} Definition 9.3. We simplify the construction through systematic use of
$\dagger$-spaces. 

We then construct a natural map from syntomic cohomology to locally analytic de Rham cohomology. This is used to show that syntomic Chern classes (and hence \'etale Chern classes) factor through locally analytic group cohomology.

\subsection{Weakly formal schemes and $\dagger$-spaces}

We review the properties of Gro\ss{}e-Kl\"onne's theory of $\dagger$-spaces that we need. See \cite{GK}, \cite{GK2} for the complete treatment. Loosely a $\dagger$-space is a rigid analytic space with structure sheaf of {\em overconvergent functions}. If $X^\dagger$ is a $\dagger$-space, we denote $\Ohdagger{X}$ the ring of overconvergent functions on $X$.

\begin{ex}\label{2.1}Let $B^\dagger$ be the ``closed unit disk'' over $K$, i.e., 
\[ B(\C_p)=\{x\in\C_p\mid |x|_p\leq  1\}.\]
 Then 
\[ 
\Ohdagger{B}=\left\{ f(t)=\sum_{n=0}^\infty a_nt^n\mid a_n\in K,\text{$f$ convergent on $|x|\leq 1+\epsilon$ for some $\epsilon$} \right\}
\]
On the other hand, let $\Delta^{\dagger}$ be the ``open unit disk'' over $K$, i.e., 
\[ \Delta^{\dagger}(\C_p)=\{x\in\C_p\mid |x|<1\}.\]
 Then
\[ 
\Ohdagger{\Delta}=\left\{ f(t)=\sum_{n=0}^\infty a_nt^n\mid a_n\in K,\text{$f$ convergent on $|x|<1$} \right\}
\]
\end{ex}

Weakly formal schemes play the role in the theory of $\dagger$-spaces which formal schemes play in
the theory of rigid analytic spaces.

Let $\Xf$ be a weakly formal $R$-scheme (\cite{GK} Kapitel 3).
We denote $\Xf^\dagger$ its {\em generic fibre} (loc. cit. Korollar 3.4) as $\dagger$-space. Let $\Xf_k$ be the special fibre of $\Xf$. This is a $k$-scheme locally of finite type. There is
a natural {\em specialisation map} 
\[ \sp: \Xf^\dagger\to \Xf_k\]
 on the underlying topological spaces.

\begin{ex}\label{schemedagger}There is a natural functor
\[ (\widehat{\cdot}):\text{$R$-schemes of finite type}\to \text{ weakly formal $R$-schemes}\]
the weak completion of the special fibre. It preserves special fibres. For an $R$ scheme of finite type, we put
\[ X^\dagger=(\widehat{X})^\dagger.\] 
Note that this is {\em not} the $\dagger$-space attached to the variety $X_K$ if $X$ is not proper.
We have $X^\dagger(K)=X(R)$.
\end{ex}

\begin{defn}A {\em $\dagger$-space with reduction} is a triple 
$\Xf=(\Xf^\dagger,\Xf_k,\sp)$ of a $\dagger$-space $\Xf^\dagger$, a $k$-scheme locally of finite type $\Xf_k$ and a continuous map $\sp:\Xf^\dagger\to \Xf_k$ of the
underlying topological spaces. Morphisms are defined in the obvious way. A morphism $\Xf\to \Yf$ of $\dagger$-spaces is called closed immersion if it is a closed immersion on both components. It is called smooth if it is smooth on both components.
\end{defn}
As discussed before, any weakly formal scheme $\Xf$ gives rise to a $\dagger$-space with reduction. 
If $\Pf$ is another weakly formal $R$-scheme,  $Y_k$ a $k$-scheme locally of finite type and
$Y_k\to \Pf_k$  a closed immersion, then the {\em tubular neighbourhood}
$]Y_k[_{\Pf}=\sp_{\Pf}^{-1}(Y_k)$ of $Y_k$ in $\Pf_K$ is also a natural $\dagger$-space with reduction. 

\begin{rem}The tubular neighbourhood $]Y_k[_\Pf$ with its reduction should be induced by some weakly formal scheme, namely the weak completion of $Y_k$ in $\Pf$. However, the theory of $\dagger$-spaces has not yet been developed up to this point.
\end{rem}

\begin{ex}For $Y=\Spec R$, the $\dagger$-space $Y^\dagger$ consists of one point.
 For $P=\A^1$, we have 
$P^\dagger=\{x\in\C_p\mid |x|\leq 1\}=\Delta^\dagger$ (see Example \ref{2.1}).
 For an $R$-valued point $a:Y\to  P$, we have
$]Y_k|_{P}(\C_p)=\{ x\in\C_p\mid |x-a|<1\}$.
\end{ex}

Smooth $\dagger$-spaces have a well-behaved theory of differential forms (\cite{GK2} 4.1). For a $\dagger$-space $Y^\dagger$, let
$\Omega^*_{Y^\dagger}$ be the complex of sheaves of  (overconvergent) differential forms on $Y^\dagger$. 
For a closed immersion of $\dagger$-spaces $Z^\dagger\to Y^\dagger$ with ideal of definition $I$, let
\[ \Fil^n_{Z^\dagger}\Omega^*_{Y^\dagger}=I^n\to I^{n-1}\Omega_{Y^\dagger}^1\to I^{n-2}\Omega_{Y^\dagger}^2\to \dots\]
be the {\em Hodge filtration}. For $Z^\dagger=Y^\dagger$ this yields the stupid filtration $\Omega^{\geq n}_{Y^\dagger}$. 
The complexes $\Fil^n_{Z^\dagger}\Omega^*_{Y^\dagger}$ are functorial with respect to such pairs.

If $\Yf\to\Pf$ is a closed immersion of weakly formal schemes, then 
$\Yf^\dagger\to ]\Yf_k[_{\Pf}$ is a closed immersion of $\dagger$-spaces.

\begin{prop}[Poincar\'e Lemma]\label{poincare}
Let $i:\Xf\to \Pf$ and $i':\Xf\to \Pf'$ be  closed immersions of weakly formal schemes with $\Pf,\Pf'$ smooth. Let $u:\Pf'\to\Pf$ be a smooth morphism compatible with the inclusion of $\Xf$, i.e., $u\circ i'=i$.
Then 
\[ \Fil^n_{\Xf^\dagger}\Omega^*_{]\Xf_k[_{\Pf}}\to
u_*\Fil^n_{\Xf^\dagger}\Omega^*_{]\Xf_k[_{\Pf'}}=Ru_*\Fil^n_{\Xf^\dagger}\Omega^*_{]\Xf_k[_{\Pf'}}
\]
 are quasi-isomorphisms of complexes of sheaves. 
\end{prop}
\begin{proof}
The cohomological assertion depends on a weak fibration formula as in rigid cohomology, \cite{Ber} 1.3.2.

It suffices to consider the case $\Xf$, $\Pf$, $\Pf'$ affine. Then 
$\Xf^\dagger$, $\Pf^\dagger$ and $(\Pf')^\dagger$ are affinoid. By making them small enough we can assume that the conormal bundle of $\Xf$ in $\Xf'=\Xf\times_\Pf\Pf'$ is free of rank $d$. Let
$t_1,\dots,t_d\in \Oh(\Pf')$ be a regular sequence defining $\Xf$ in $\Pf$. Then $dt_1,\dots,dt_d$
are a basis of $\Omega^1_{\Pf'/\Pf}$ in a neighbourhood of $\Xf$. These sections define
a morphism
\[ \Pf'\to \Pf''=\Pf\times(\hat{\A}^1)^d\]
\'etale in a neighbourhood of $\Xf$. The closed immersion of $\Xf$ is given by the zero-section on
$(\hat{\A}^1)^d$.

By  \cite{Ber} Proposition 1.3.1
\[ ]\Xf_k[_{\Pf'}\to ]\Xf_k[_{\Pf''}\]
is an isomorphism for the corresponding rigid analytic varieties. By \cite{GK2} Theorem 1.12 (a) this implies that the morphism of dagger-spaces is an isomorphism. On the other hand
\[ 
]\Xf_k[_{\Pf''}\isom]\Xf_k[_{\Pf}\times (\A^{1\dagger})^d= ]\Xf_k[_{\Pf}\times\Delta^d
\]
where $\Delta$ is the open unit disc. (There is an obvious map to the right hand side which is an isomorphism of rigid analytic spaces and hence also of dagger-spaces). 

We now turn to the statement on differential forms. 
By \cite{GK2} Satz 4.12 
\[ \Omega^*(]\Xf_k[_{\Pf}\times\Delta^d)=\Omega^*(]\Xf_k[_{\Pf})\tensor\Omega^*(\Delta)^{\tensor d}\]
The filtration is compatible with this
decomposition. This reduces the proof to the case $\Pf'=\Delta$ and $\Xf=\Pf$ the zero-section.

As $\Delta$ is Stein, $Ru_*=u_*$.
Let $t$ be the parameter of $\Delta$. The filtration has two steps:
\[ \Fil^0\Omega^*(\Delta)=[\Oh(\Delta)\to\Omega^1(\Delta)]\hspace{3em}  
\Fil^1\Omega^*(\Delta)=[t\Oh(\Delta)\to\Omega^1(\Delta)]
\]
The differential is an isomorphism on $\Fil^1$, i.e. the complex is acyclic. The kernel on $\Fil^0$ 
consists of constants functions, i.e. the cohomology of a single point.
\end{proof}

\begin{prop}[Rigid Cohomology ]\label{GK}
Let $\Xf\to \Pf$ be a closed immersion of smooth weakly formal schemes. Then
$H^i( ]\Xf_k[_{\Pf}, \Omega^*)$ is naturally isomorphic to rigid cohomology of $\Xf_k$ in the sense of
Berthelot.
\end{prop}
\begin{proof}This is \cite{GK} Proposition 8.1 (b) or \cite{GK2} Theorem 5.1.\end{proof}

\subsection{Syntomic cohomology}
We define syntomic cohomology on affine $\dagger$-spaces with reduction. The case of most interest is the one of the weak completion of an affine $R$-scheme.
The restriction to the affine case is not essential. It simplifies the construction slightly because all $\dagger$-spaces which occur are acyclic for cohomology of coherent sheaves.

\begin{defn}\label{syndef}A {\em syntomic data} for an affine $\dagger$-space with reduction $\Xf$ is a collection of
\begin{itemize}
\item a smooth  affine weakly formal $R_0$-scheme $\Pf_0$ together with a $\sigma$-linear lift $\Phi$ of
absolute Frobenius on the special fibre $\Pf_{0k}$;
\item a closed immersion $\Xf_k\to \Pf_{0k}$;
\item a smooth affine weakly formal $R$-scheme $\Pf$;
\item a closed immersion $\Xf\to \Pf$ of $\dagger$-spaces with reduction and a morphism of weakly formal schemes $\Pf_0\to \Pf$ such that the two
maps $\Xf_k\to \Pf_k$ agree.
\end{itemize}
A morphism of syntomic data for $\Xf$ is a pair of smooth morphisms 
$u_0:\Pf'_0\to \Pf_0$, $u:\Pf'\to \Pf$ such that the obvious diagrams commute.

Let $n\in\Z$.
The {\em syntomic complex} $R\Gamma_\syn(\Xf,n)_{\Pf_0,\Pf}$ attached to this data is defined as
\[ 
  \cone\left( \Fil^n_{\Xf^\dagger}\Omega^*(]\Xf_k[_{\Pf})\oplus \Omega^*(]\Xf_k[_{\Pf_0})\to
    \Omega^*(]\Xf_k[_{\Pf})\oplus \Omega^*(]\Xf_k[_{\Pf_0})\right)[-1]
\]
where the map is given by $(a,b)\mapsto (a-b,(1-\Phi^*/p^n)b)$. Its cohomology
is called {\em syntomic cohomology} $H^i_\syn(\Xf,n)$ of $X$.
If $X$ is a smooth affine scheme, $H^i_\syn(X,n)$ is defined as defined as syntomic cohomology of its weakly formal completion.
\end{defn}

\begin{rem}
By the Poincar\'e Lemma \ref{poincare}, a morphism of syntomic data induces an isomorphism on syntomic cohomology, i.e., a quasi-isomorphism of syntomic complexes. Note, however, that
the system of syntomic data is not filtering:  a pair of morphism of syntomic data
$\alpha,\beta:(\Pf_0,\Phi,\Pf)\to (\Pf'_0,\Phi',\Pf')$ is not equalised on a third data. To obtain
a complex independent of choices (and hence a functorial theory), one has to proceed as Besser
in \cite{Be} Definition 4.11 - Definition 4.13. We do not go into the details.
\end{rem}

\begin{rem}The restriction to the affine case is not necessary. In the general
case one has to replace the global sections $\Omega^*(]\Xf_k[_\Pf)$ by global sections of a functorial injective resolution of $\Omega^*_{]\Xf_k[_\Pf}$.
\end{rem}

\begin{ex}\label{pointcase}
Let $\Xf=\Sp \widehat{R}$ the weakly formal completion of $\Spec R$. Then $\Pf_0=\Sp \widehat{R}_0$ with $\Phi=\sigma$ and $\Xf=\Pf$ is a syntomic data.
We have $]\Xf_k[_{\Xf}=\Xf^\dagger$ (a single point) and hence the ideal of definition $I$ vanishes. For $n>0$,
the complex $\Omega^{\geq n}(\Xf^\dagger)$ vanishes. Moreover, $\Omega^0(\Xf^\dagger)=K$
(constant functions on a single point).
Hence the syntomic complex $R\Gamma_\syn(\Spec R,n)_{\Pf_0,\Xf}$ is simply
\[ \cone\left(K_0\xrightarrow{(1,1-\sigma/p^n)}K\oplus K_0\right)[-1]\]
The $\Q_p$-linear map $(1-\sigma/p^n)$ is bijective, hence 
\[ \eta^{-1}: K[-1]\to R\Gamma_\syn(\Spec R,n)_{\Pf_0,\Xf}\]
is a quasi-isomorphism. Hence for $n>0$
\[
  H_\syn^i(\Spec R,n)=\begin{cases}K&i=1\\
                              0&\text{otherwise}
                 \end{cases}
\]
\end{ex}
This identification will be used very often in the sequel.
\begin{defn}\label{eta}Let $n>0$. We denote
\[ \eta:H^1(R,n)\to K\] 
the isomorphism of Example \ref{pointcase}.
\end{defn}

\begin{rem}Let $X$ be a smooth affine $R$-scheme, $a\in X(R)$, $c\in H^1_\syn(X,n)$. Then
$\eta(a^*c)\in K$. This means that $c$ induces a (set-theoretic) map
\[ c:X(R)\to K\]
We are going to show in the next section that this map is in fact locally analytic on $X(R)$
(but not necessarily rigid analytic in general).
\end{rem}

\begin{prop}\label{Besser}For $R$-schemes $X$, syntomic cohomology as defined above agrees with $\tilde{H}^i_\syn(X,n)$ defined
by Besser \cite{Be} Definition 9.3. For $R=R_0$, it agrees with syntomic cohomology defined
by Gros \cite{G}.
\end{prop}
\bew Besser uses direct limits over rigid analytic functions in strict neighbourhoods of
$]\Xf_k[_{\Pf}$ rather than $\dagger$-spaces.
By Proposition \ref{GK} this amounts to the same. Apart from this point, the definitions agree.
The second statement is \cite{Be} Proposition 9.4.
\bewende

\begin{rem}
The theory immediately extends to simplicial schemes over $R$, in particular to $B\GL_N$.
\end{rem}

\begin{prop}[\cite{Be} 9.10, \cite{Ni2}]\label{synet} Let $X$ be a smooth affine $R$-scheme. Then there is a natural morphism
\[ H^i_\syn(X,n)\to H^i_\et(X_k,\Q_p(n))\]
of cohomology theories compatible with Chern classes.
\end{prop}
\bew Besser constructs such a transformation for his version of syntomic cohomology.
(The main step in the proof is due to from Niziol, see  \cite{Ni2}.) 
In his proof the map factors by construction through his $\tilde{H}^i_\syn(X,n)$ which is
our $H^i_\syn(X,n)$ (Proposition \ref{Besser}). 
\bewende

\subsection{Evaluation maps and Chern classes}

Let $X$ be a smooth affine $R$-scheme.
\begin{defn}Let 
\[
 X^\delta=\coprod_{a\in X(R)_\cdot}\Spec R
\]
\end{defn}
Recall that $\Ohdeltatwo{X}{F}$ denotes $F$-valued set-theoretic functions on
$X(R)$.
\begin{lemma}\label{eval}Let $n>0$. The natural morphism of schemes $X^\delta\to X$ induces
a natural map
\[ \eval:R\Gamma_\syn(X,n)\to \Ohdeltatwo{X}{H^1_\syn(\Spec R,n)}[-1]\stackrel{\eta}{\to} 
  \Ohdelta{X}[-1]\]
with $\eta$ as in Definition \ref{eta}.
\end{lemma}
\bew
We first use the functoriality of syntomic complexes
\[ R\Gamma_\syn(X,n)\to R\Gamma_\syn(X^\delta,n)\]
This gives the formula of the Lemma by Example \ref{pointcase}.
\bewende

Applying this map to the simplicial scheme $\BG$ for our smooth affine group-scheme $G=\GL_N$, we get by definition of group cohomology a natural map
\[ \eval: H^{2n}_\syn(\BG,n)\to H^{2n-1}(G(R),K)\ .\]

Gros (for $R=R_0$) and Besser (general case) have established the existence of Chern classes in syntomic cohomology.
The key ingredient is a universal Chern class for $i\leq N$
\[ c_i\in H^{2i}_\syn(\BGLN,i)\]
It is uniquely characterised by the fact that its image in $\Fil^iH^{2i}(B\GL_{N,K}^\dagger,\Omega^*)$
is the usual Chern class in de Rham cohomology (\cite{Be} Proposition 7.4 and the discussion following it). 

\begin{defn}[\cite{Be} Theorem 7.5]\label{syntomicchern}The syntomic Chern class 
\[ c_n\in H^{2n-1}(\GL_N(R),H^1_\syn(R,n))\]
is given by applying the evaluation map of Lemma \ref{eval} for $\BGLN$ to
the universal Chern class.
\end{defn}

We denote $R\Gamma_\et(K,\Q_p(n))$ the complex computing continuous \'etale cohomology of $K$ with coefficients in $\Q_p(n)$. This agrees with continuous cohomology of the group $\Gal(\bar{K}/K)$.
\begin{prop}\label{etale}
As before let $G=\GL_N$ considered as smooth group scheme over $R$. For $n>1$,  there is a natural commutative diagram
\[\begin{xy}
\xymatrix{
H^{2n}_\syn(\BG,n)\ar[r]\ar[d]&H^{2n}_\et(\BG_K,\Q_p(n))\ar[d]\\
H^{2n-1}(G(R),H^1_\syn(R,n))\ar[r]\ar[rd]_\eta& H^{2n-1}(G(R),H^1_\et(K,\Q_p(n)))\\
&H^{2n-1}(G(R),K)\ar[u]_{\exp_\BK}
}
\end{xy}\]

For $n=1$ the diagram reads

\[\begin{xy}
\xymatrix{
H^{2}_\syn(\BG,1)\ar[rr]\ar[d]&&H^{2}_\et(\BG_K,\Q_p(1))\ar[d]\\
H^{1}(G(R),H^1_\syn(R,1))\ar[r]\ar[rd]_\eta& H^{1}(G(R),H^1_\et(K,\Q_p(1)))\ar[r]& H^{1}(G(R),R\Gamma_\et(K,\Q_p(1)))\\
&H^{1}(G(R),K)\ar[u]_{\exp_\BK}
}
\end{xy}\]

For $1\leq n$ and $N$ big enough, the universal syntomic Chern class is mapped
to the universal \'etale Chern class.
\end{prop} 
\bew The vertical maps are the ones from Lemma \ref{eval} and their \'etale
analogue respectively. Note that $R\Gamma_\et(K,\Q_p(n))$ is concentrated in
degrees $1,2$ for $n\neq 0$ (even in degree $1$ for $n\neq 0,1$). 
The natural transformation of Proposition \ref{synet} gives the horizontal maps.
This yields the upper commutative square. By loc. cit. it is compatible with Chern classes. 

In \cite{Be} Proposition 9.11
the relation to the exponential is made explicit. This gives the lower triangle.
\bewende

\begin{rem} This proposition reduces the proof of our Main Theorem~\ref{mainthm} to a statement on universal Chern classes in syntomic cohomology.
\end{rem}

\subsection{Analyticity of evaluation maps}
We are going to show that the syntomic Chern classes $c_n$ is an element of the
continuous group cohomology of $\GL_N(R)$. We want to prove:
\begin{thm}\label{evalthm}Let $X$ be a smooth affine $R$-scheme, $n\geq 1$. Then
the evaluation map  of Lemma \ref{eval} factors naturally via
a morphism in the derived category of $\Q_p$-vector spaces
\[ \eta:R\Gamma_\syn(X,n)\to \Omega^{<n}(X^\la)[-1]\to \Ohla{X}[-1]\]
Moreover, $\eta$ is represented by a natural sequence of morphisms of  complexes and formal inverses of quasi-isomorphisms of complexes. For $X=\Spec R$, the map agrees with the one defined previously (see Definition \ref{eta}).
\end{thm}

\begin{defn}Let $Y$ be a smooth affine $R$-scheme, 
Let $\mathring{Y}$ the $\dagger$-space with reduction with special and 
generic fibre 
\begin{align*}
\mathring{Y}_k&=\coprod_{a\in Y(k)}\Spec k\\
 \mathring{Y}^\dagger&=\coprod_{a\in Y(k)}]a[_{Y}
\end{align*}
together with the natural specialisation map.
\end{defn}

\begin{rem}Note that by definition $(\mathring{Y})^\dagger(K)=Y^\dagger(K)=Y(R)$. 
The locally analytic manifolds on these spaces agree.
\end{rem}

\begin{lemma}\label{ybeta}Let $(\Pf_0,\Phi,\Pf)$ a syntomic data for $Y$. Then
$(\Pf'_0,\Phi',\Pf')$ with
$\Pf'=  \coprod_{a\in Y(k)}\Pf_0$, $\Pf'=\coprod_{a\in Y(k)}\Pf$
with $\Phi'$ operating as $\sigma$ on $Y(k)$ and via $\Phi$ on the $\Pf_0$
is
a syntomic data for $\mathring{Y}$.
Let $n>0$. Then there is a natural isomorphism
\[ R\Gamma_\syn(\mathring{Y},n)_{(\Pf'_0,\Pf')}\to\Omega^{< n}(\mathring{Y}^\dagger)[-1]\]
in the derived category of $\Q_p$-vector spaces. Moreover, the morphism is represented by a natural
sequence of quasi-isomorphisms of complexes going either direction.
\end{lemma}
\bew 
 All  properties of a syntomic data follow from functoriality. 

By the Poincar\'e Lemma \ref{poincare}, the natural inclusion
\[ \bigoplus_{a\in Y(k)}K_0\to \Omega^\dagger(]\mathring{Y}_k[_{\Pf'_0})\]
is a quasi-isomorphism. Hence
\[ \cone\left(\Fil^n_{\mathring{Y}^\dagger}\Omega^*(]\mathring{Y}_k[_{\Pf'})\oplus  \bigoplus_{a\in Y(k)}K_0
 \to \Omega^*(]\mathring{Y}_k[_{\Pf'})\oplus  \bigoplus_{a\in Y(k)}K_0\right)[-1]
\]
is quasi-isomorphic to $R\Gamma_\syn(\mathring{Y},n)$.
The map $\Phi^*$ operates as $\sigma$ on $K_0$ and permutes the elements of $Y(k)$. This map
preserves the norm on $\bigoplus_{a\in Y(k)}K_0$. Hence
 the map $1-\Phi^*/p^n$ is 
an isomorphism for $n>0$. Hence the inclusion of
\[ \cone\left(\Fil^n_{\mathring{Y}}\Omega^*(]\mathring{Y}_k[_{\Pf'})\to \Omega^*(]\mathring{Y}_k[_{\Pf'})\right)[-1]
\]
is a quasi-isomorphism. Now we apply the Poincar\'e Lemma again: 
$\Fil^n_{\mathring{Y}}\Omega^*(]\mathring{Y}_k[_{\Pf'})$ is quasi-isomorphic to the subcomplex
\[ \Fil^n_{\mathring{Y}^\dagger}\Omega^\dagger(\mathring{Y}^\dagger)=\Omega^{\geq n}(\mathring{Y}^\dagger) \]
and $\Omega^\dagger(]\mathring{Y}_k[_{\Pf'})$ to $\Omega^\dagger(\mathring{Y}^\dagger)$.
Finally the cone 
\[ \cone\left(\Omega^{\geq n}(\mathring{Y}^\dagger)\to\Omega^\dagger(\mathring{Y}^\dagger)\right)[-1]\]
 is  quasi-isomorphic to the quotient.
\bewende

\begin{proof}[Proof of Theorem \ref{evalthm}]
We define $\eta$ as the composition of the natural map
\[ R\Gamma_\syn(Y,n)\to R\Gamma_\syn(\mathring{Y},n)\]
with the quasi-isomorphism of Lemma \ref{ybeta}
\[ R\Gamma_\syn(\mathring{Y},n)\to \Omega^{< n}(\mathring{Y}^\dagger)[-1]\]
Note finally that $\mathring{Y}^\dagger(K)=Y(R)$ as $K$-manifolds and that
overconvergent differentials are locally analytic.

Now assume $Y=\Spec R$. Then $\mathring{Y}=Y$, $\Oh(\mathring{Y}^\dagger)=K$. The
the chain of isomorphisms in Lemma \ref{ybeta} agrees with the one in
Example \ref{pointcase}.

As $\eta$ is natural, this implies that the evaluation map of Lemma \ref{eval} factors through
$\eta$. 
\end{proof}

\begin{comment}

\begin{proof}[Proof of Theorem \ref{evalthm}]
Let $X_\cdot$ be a smooth affine simplicial scheme, $n>0$.
By construction, $\eta$ is not only a morphism in the derived category but represented
by explicit natural morphisms and inverses of quasi-isomorphisms. Hence it extends to
simplicial schemes. We get
\[ \eta:R\Gamma_\syn(X_\cdot,n)\to \Oh(\mathring{X}^\dagger_{K\cdot})[-1]\]
For any smooth $Y$, we have natural inclusions
\[ \Oh(Y^\dagger)\to \Oh_\la(Y(R)^\an)\to \maps(Y(R),K)\]
where $\Oh_\la$ denotes locally analytic $K$-valued functions on a $p$-adic manifold and
$\maps$ denotes set-theoretic maps.
Note that $H^i_\Delta(X(R)_\cdot,K)$ is cohomology of $\maps(X(R)_\cdot,K)$. By naturality
of $\eta$, the evaluation map factors through $\eta$ and hence via
cohomology of $\Oh_\la(X_\cdot(R)^\an)$ as claimed.
\bewende
\end{comment}

We apply the arguments to $\BG$. 
For later use we record a couple of commutative diagrams:
\begin{prop}\label{diagrams}
As before let $G=\GL_N$ as smooth affine algebraic group over $R$. 
Let $\cG=1+\pi M_N(R)$ as locally analytic $K$-manifold. We have $\cG=\cG^\dagger(R)$ where
 $\cG^\dagger=]e[_G\subset G^\dagger$ as dagger-space. 
Then the following diagram commutes:
\[\begin{xy}
\xymatrix{
H^{2n}(\Omega^{\geq n}(\BG^\dagger))\ar[r]&H^{2n}(\Omega^{\geq n}(\BcG^\dagger))\\
H^{2n}_\syn(\BG,n)\ar[u]^d\ar[d]^\eval\ar[r]^\eta & H^{2n-1}(\Omega^{<n}(\BcG^\dagger))\ar[u]_\partial\ar[d]\\
H^{2n-1}(\cG,K) & H^{2n-1}_\la(\cG,K)\ar[l]
}
\end{xy}\]
where $d$ is induced from the natural map 
$R\Gamma_\syn(\BG,n)\to \Omega^{\geq n}(\BG^\dagger)$ (see Definition \ref{syndef}),  $\partial$ is induced from the connecting map of the short exact sequence of complexes 
\[ 0\to \Omega^{\geq n}(\BcG^\dagger)\to \Omega^*(\BcG^\dagger)\to \Omega^{<n}(\BcG^\dagger)
\to 0\] 
$\eval$ is the evaluation map \ref{eval} and $\eta$ is the map of Theorem \ref{evalthm}.

For $1\leq n$ and $N$ big enough, we define the Chern class $c_n^\an\in H^{2n}(\Omega^{\geq n}B.\cG^\dagger)$ 
as image of the Chern class in algebraic de Rham cohomology. 
Then the  universal syntomic Chern class is mapped to $c_n^\an$ in the top right corner.
\end{prop}
\bew 
We apply Lemma \ref{ybeta} to $\BG$ and restrict to $\BcG\subset B.\mathring{G}$. This
defines the map $\eta$ in the middle.
The commutativity of the lower square follows from Theorem \ref{evalthm} together
with the definition of locally analytic group cohomology.

By definition of syntomic Chern classes they are mapped to the standard
Chern classes in algebraic and hence also overconvergent de Rham cohomology, see \cite{Be} Proposition 7.4 and the discussion following it.
(Note that Besser uses a more refined version of syntomic cohomology than we do. His 
version of $d$ takes values is {\em algebraic} de Rham cohomology.)

It remains to check commutativity of the upper square. By construction of $\eta$, we have (for $Y$ a smooth affine $R$-scheme):
\[\begin{xy}
\xymatrix{
\Omega^{\geq n}(Y^\dagger)\ar[r]& \Omega^{\geq n}(\mathring{Y}^\dagger)\\
R\Gamma_\syn(Y,n)\ar[u]^d\ar[r]& R\Gamma_\syn(\mathring{Y})\ar[u]^d\ar[r]& 
\cone\left(\Omega^{\geq n}(\mathring{Y}^\dagger)\to \Omega^*(\mathring{Y}^\dagger)\right)[-1]\ar[ul]_\partial
}
\end{xy}\]
Restriction to $\BcG$ gives the statement.
\end{proof}
\begin{remark} 
$\cG$ is an open and closed 
subgroup of  finite index in $\GL_N(R)$, hence group cohomology
of $\cG$ and $G(R)$ agree (with rational coefficients).
\end{remark}

\begin{remark} The above diagrams work without changes for all smooth algebraic group schemes over
$R$. 
\end{remark}

\subsection{Analyticity of Chern classes}
\begin{thm}\label{syntomicfactors}Let $n>0$ and $N>n$. There exists an element
\[ (\eta c_n)^\la\in H^{2n-1}_\la(\GL_N(R),K)\]
which has the same image in $H^{2n-1}(\GL_N(R),K)$ as the syntomic Chern classes
\[ c_n\in H^{2n-1}(\GL_N(R),H^{1}_\syn(R,n)) \]
(cf. Definition \ref{syntomicchern})
under the map $\eta:H^1_\syn(R,n)\to K$. In particular, $c_n$ is the image of an
element
\[ c_n^\cont\in H^{2n-1}_\cont(\GL_N(R),H^{1}_\syn(R,n))\]
\end{thm}
\begin{proof} 
Apply Proposition \ref{diagrams} to $\BG$.
The image of  the universal Chern class $c_j$ in $H^{2j-1}(\GL_N(R),H^1_\syn(R,j)$
is by Definition \ref{syntomicchern} the syntomic Chern class for $R$. By the diagram in Proposition \ref{diagrams}
it is the image of a locally analytic class.
In particular it is the image of a continuous class.
\end{proof}

\begin{ex}Let $N=1$, $j=1$. 
By \cite{G} Proposition 4.1, the first Chern class
\[ c_1:R^*\tensor \Q_p\to H^1_\syn(R,1)=K\]
is given by the $p$-adic logarithm $\log_p$. This function is locally analytic but not rigid analytic. 
However, it is (overconvergent) rigid analytic on the open unit disc $\Delta^\dagger$.
Our proof shows the same behaviour also for higher Chern classes.
\end{ex}
 
\begin{cor}\label{chernanalytic}Let $n>0$. The \'etale Chern class 
\[ c_n\in H^{2n-1}(\GL_N(R),H^1_\et(K,n))\]
is the image of an element
\[ c_n^\cont\in H^{2n-1}_\cont(\GL_N(R),H^1_\et(K,n))\]
in continuous group cohomology.
\end{cor}
\bew Combine Theorem \ref{syntomicfactors} with Proposition \ref{etale}.
\bewende

Theorem \ref{syntomicfactors} allows to reduce the proof of our Main Theorem \ref{mainthm} to continuous group cohomology.

\section{The suspension map and locally analytic group cohomology}
As before let $G=\GL_N$ as algebraic group over $R$, le $\cG^\dagger=]e[_G$ the residue disc of $e$ in 
$G^\dagger$ viewed as dagger-space,
and  $\cG=\cG^\dagger(K)=1+\pi M_N(R)$ as 
$K$-Lie group. Let $\g=\gl_N$ be the $K$-Lie algebra of $G$ (see \ref{lie-algebra}).

In the last section, we constructed an element 
\[ (\eta c_n)^\la\in H^{2n-1}_\la(\GL_N(R),K)=H^{2n-1}_\la(\cG,K)\]
(see Theorem \ref{syntomicfactors}). In this section we are going to define (see Definition \ref{psidefn}) a natural map
\[ \Psi: H^{2n-1}_\la(\cG,K)\to H^{n-1}(\g,K)\ .\]
Eventually, we want to show that $\Psi( (\eta c_n)^\la )=p_n$, the primitive element in Lie algebra cohomology. 

The aim of this section is to embed $\Psi$ into a huge commutative diagram relating it 
 to the suspension for $\BG$. As the image of the universal Chern class in algebraic de Rham cohomology under this suspension map is precisely $p_n$, this will allow to deduce the claim (see Section \ref{proofsection}).

We follow
closely the ideas of Beilinson \cite{Be} as outlined by Rapoport \cite{Ra} and Burgos \cite{Burgos}.

\subsection{A commutative diagram}
We state the result of this chapter. 

\begin{thm} \label{suspdiagram}There is a  natural commutative diagram
\[
\begin{xy}\xymatrix{
	H^{2n}(\Omega^{\ge n}(B.G^\alg))\ar[rd]^\inf\ar[d]\ar^{s_G}[rr]&& 
	         H^{2n-1}_\DR(G^\alg)\ar^{\rho}_\isom[ddd]\\
	H^{2n}(\Omega^{\ge n}(\BcG^\dagger))\ar[r]^\inf &H^{2n}(W^{\ge n,\cdot}(\g))\ar[rdd]^{s_\g}\\
	H^{2n-1}(\Omega^{< n}(\BcG^\dagger))\ar[r]^\inf\ar[u]\ar[d] &H^{2n-1}(W^{< n,\cdot}(\g))\ar[u]\ar[dr]\\
	H^{2n-1}_\la(\cG,K)\ar[rr]^\Psi&&H^{2n-1}(\g,K).}
\end{xy}
\]
\end{thm}
The suspension map $s_G$ will be introduced in Section \ref{suspdefn}, $\rho$ in Lemma \ref{rhodefn} and
the algebraic Lazard isomorphism $\Psi$ in Definition \ref{psidefn}.
 The Weil algebra $W^{*,\cdot}(\g)$ and
the map $s_\g$ are defined in Section \ref{weil}. The various maps $\inf$ will be introduced in Section \ref{infdefn}.
 Finally the proof of the Theorem will given in Section \ref{proofofsuspensiondiag}.
\begin{rem} The same arguments yield the above diagram in the case of a
reductive group over $K$.
\end{rem}
\subsection{The suspension map $s_G$}\label{suspdefn}

Consider the simplicial schemes
$E.G$ and $B.G$ over $R$. 
Let $G.$ be the constant simplicial scheme and $\Delta:G.\to E.G$ be the diagonal 
inclusion. Then we have a fibre diagram
\[
	\begin{CD}
		G.@>\Delta>> E.G\\
		@.@VVV\\
		@. B.G
	\end{CD}
\]
As $E.G$ is contractible the suspension for this $G$-torsor gives us a morphism for $n>0$
\[
	s_G:H^{2n}_\DR(B.G^\alg)\to H^{2n-1}_\DR(G^\alg)
\]
(compare \cite{Burgos} Example 4.16 and recall that $H^{2n}_\DR(B.G^\alg)$ and
$H^{2n-1}_\DR(G^\alg)$ are the de Rham cohomology of the generic fibre). 
We will use another description of the suspension map in terms of the Eilenberg-Moore
spectral sequence
\[
	E_1^{p,q}=H^q_\DR(B_pG^\alg)\Longrightarrow H^{p+q}_\DR(B.G^\alg).
\]
As $E_1^{0,q}=0$ for $q>0$ we get an edge morphism
\begin{equation}\label{sGdefn}
	s_G:H^{2n}_\DR(B.G^\alg)\to E_1^{1,2n-1}=H^{2n-1}_\DR(G^\alg),
\end{equation}
which is none other than the suspension $s_G$.
In particular, the suspension is compatible with the "Hodge filtration" (defined
by $\Omega^{\ge n}(B.G^\alg)$ and $\Omega^{\ge n}(G^\alg)$) on both sides
and we get a map
\begin{equation}\label{sGfildefn}
	s_G:H^{2n}(\Omega^{\ge n}(B.G^\alg))\to H^{2n-1}(\Omega^{\ge n}(G^\alg)).
\end{equation}

\subsection{Lie algebra and de Rham cohomology}

Let $\g^\vee:=\Hom_K(\g,K)$ be the dual of $\g$.
Let $\cC^*(\g)$  be the standard complex of Lie algebra cohomology with coefficients in $K$ and
\[
	H^i(\g):=H^i(\g,K).
\]
Let $\Omega^*(G^\alg)$ be the de Rham complex of the generic fibre $G\times_RK$ of
$G$. Identifying $\frg$ with the left invariant vector fields on $G$, one has an embedding
\begin{equation}\label{lieincl}
	\g^\vee\subset \Omega^1(G^\alg)
\end{equation}
of the dual of $\g$ into the $1$-forms on $G\times_RK$, which induces a
map of complexes
\[
	\cC^*(\g)\subset \Omega^*(G^\alg).
\]
It has a splitting by evaluation at $e$.
\begin{lemma}[\cite{Ho} Lemma 4.1]\label{rhodefn}
	The above inclusion induces an isomorphism
	\[
		\rho: H^*(\g)\isom H^*_\DR(G^\alg).
	\]
\end{lemma}

\subsection{The infinitesimal version of $B.G$ and Lie algebra cohomology}
Recall that  $\Ohalg{B.G}$ is the cosimplicial ring of $K$-valued algebraic
functions on  $B.G$. If we extend
$\Spec$ in the obvious way to cosimplicial rings, we have
\[
	\Spec\Ohalg{B.G}=B.G\times_RK.
\]
We need an infinitesimal version of $B.G$. For this let $\frm_e^\alg\subset \Ohalg{B_1G}=\Ohalg{G}$ 
be the kernel of the augmentation map defined by the unit element $e\in G(R)$. 
Let $J^\cdot\subset \Ohalg{B.G}$ be the cosimplicial ideal generated by $(\frm_e^\alg)^2$. 
\begin{definition}
	Let $\Ohalg{\BGinf}:=\Ohalg{B.G}/J^\cdot$. 
	The \emph{infinitesimal version of $B.G$} is the simplicial scheme
	\[
		\BGinf:= \Spec \Ohalg{\BGinf}.
	\]
\end{definition}
The closed immersion $\BGinf\subset B.G$ is defined by the canonical map
\[
	\Ohalg{\BG}\to \Ohalg{\BGinf}.
\]
It is central for our arguments that this 
 has a counterpart for the (overconvergent) rigid analytic functions.
 \begin{lemma}\label{inf1}
The  natural map of cosimplicial rings
\[
	\inf:\Ohalg{\BG}\to \Ohalg{\BGinf}
\]
factors naturally via
\begin{equation}
	\inf:\Ohdagger{\BcG}\to \Ohalg{\BGinf}
\end{equation}
\end{lemma}
\bew
First consider more generally a smooth $R$-scheme $X$ and $x\in X(R)$. 
Let $\cX^\dagger$ an open neighbourhood of $x$ in $X^\dagger$ .
 Let $J$ be an ideal of $\Oh_{X^\alg,x}$ and put $J^\dagger=J\Oh_{\cX^\dagger,x}$.
For $J=\frm_x$ the maximal ideal, $\frm_x^\dagger$ is indeed the maximal ideal of $\Oh_{\cX^\dagger,x}$.
Hence for $J$ containing a power of $\frm_x$, 
\[ \Oh_{\cX^\dagger,x}/J^\dagger\isom \Oh_{X^\alg,x}/J\]
(Note that stalks of the structure sheaf of $\cX^\dagger$ as dagger-space and the corresponding rigid analytic variety agree.) 
All components of the simplicial scheme $\BGinf$ are of this form.
By definition 
$\BcG=]e[_{\BG}$ is an open neighbourhood of $e$ in $\BG^\dagger$. This yields the claim.
\bewende
To formulate the next two propositions, we need the concept of normalisation for
cosimplicial rings.
\begin{definition}
	For any cosimplicial object $A^\cdot$ in an abelian category
	let $\cC A^\cdot$ be the complex with $\cC A^n=A^n$ and 
	differential $d=\sum_{i=0}^{n+1}(-1)^i\delta_i$. The \emph{normalisation} is
	the subcomplex 
	\[
		\cN A^\cdot=\bigcap_{i=0}^{n-1}\ker\sigma^i.
	\]
\end{definition}
\begin{proposition}\label{normofAG}
	There is a natural isomorphism of complexes
	\[
		\cN \Ohalg{\BGinf}\isom \cC^\cdot(\g).
	\]
\end{proposition}
\begin{proof}
	This is \cite{Ra} Lemma 3.1 or \cite{Burgos} Proposition 8.9. The first reference works
	over $\bbR$ and the second over $\bbC$. But by inspection both proofs work without 
	any changes over an arbitrary field of characteristic $0$.
\end{proof}
\begin{defn}\label{psidefn}Let 
\[\Psi:\Ohdagger{\BcG}\to \cC(\g)\]
be given by the composition
$\pi^*:\Ohdagger{B_n\cG}\to\Ohdagger{E_n\cG}\isom \Ohdagger{\cG}^{\hat{\tensor}(n+1)}$
with the map 
\begin{align*}
\Ohdagger{\cG}^{\hat{\tensor}(n+1)}&\to \bigwedge^n\g^\vee\\
f_0\tensor\dots\tensor f_n&\mapsto f_0(e)df_1(e)\wedge\dots\wedge df_n(e)
\end{align*}
Here $\pi:E_n\cG^\dagger\to B_n\cG^\dagger$ is the map
$(g_0,\dots,g_n)\mapsto (g_0g_1^{-1},\dots,g_{n-1}g_n^{-1})$ 

We also denote $\Psi$ the analogous maps on $\Ohalg{\BG}$ and $\Ohla{\BG}$.
The induced map on cohomology is called {\em algebraic Lazard morphism}.
\end{defn}

\begin{lemma}\label{algebraicLazard}
	The morphism $\Psi$ agrees with the natural morphism of complexes 
	\[
		\cN\Ohdagger{B.G}\to\cN \Ohalg{\BGinf}\isom \cC^\cdot(\g)
	\]
	induced by $\inf$, where the isomorphism is the one from Proposition \ref{normofAG}.
\end{lemma}
\bew 
Burgos (cf. \cite{Burgos} Theorem 8.4.) first defines a map 
\begin{align}
	\Ohalg{E_nG}&\to \Omega^n_{G/K}\\
	f_0\otimes\dots\otimes f_n&\mapsto f_0df_1\wedge\dots\wedge df_n,
\end{align}
which factors over $\Ohalg{\EGinf}$.
As $\cC^n(\frg)\subset \Omega^n_{G/K}$
are the $G$-invariant differential forms it factors naturally as
\[
\begin{CD}
	\Ohalg{E_nG}@>>>\Omega^n_{G/K}\\
	@A\pi^*AA@AAA\\
	\Ohalg{B_nG}@>\Psi>>\cC^n(\g)\\
\end{CD}
\]
Evaluation at $e$ is a splitting of the right vertical map. 
As $\Psi$ factors through $\Ohdagger{B_nG}$ and $\Ohalg{\EGinf}^G=\Ohalg{\BGinf}$,
we get the desired result.
 \bewende

\subsection{The Weil algebra and the infinitesimal suspension map $s_\g$}\label{weil}

In this section we will define the bigraded Weil algebra $W^{*,\cdot}(\g)$ and 
the suspension on the level of the Lie algebra $\g$.

\begin{definition}\label{weildefn}
	The {\em Weil algebra} is the bigraded algebra with
	\[
		W^{p,q}(\g):=\Sym^p\g^\vee\otimes\Lambda^{q-p}\g^\vee
	\]
	and $W^{p,q}(\g)$ has total degree $p+q$ (this means that $\Sym^p\g^\vee$ has degree $2p$).
	Write $W^n(\g):=\bigoplus_{p+q=n}W^{p,q}(\g)$.
\end{definition}
The Weil algebra has a differential $\delta:W^n(\g)\to W^{n+1}(\g)$,
which is uniquely determined on $\Sym^1\g^\vee$ and $\Lambda^{1}\g^\vee$ as follows:
let $X_1,\dots, X_k$ be a basis of $\g$ and $X^\vee_1,\dots, X^\vee_k$ be the dual basis and
$h:\Sym^1\g^\vee\to \Lambda^{1}\g^\vee$ the identity map, then
\begin{align}
	\delta:\Lambda^{1}\g^\vee&\to \Sym^1\g^\vee\oplus \Lambda^{2}\g^\vee\\
	\nonumber	X^\vee&\mapsto h(X^\vee)+dX^\vee,
\end{align}
where $d:\cC^\cdot(\g)\to \cC^{\cdot+1}(\g)$ is the differential in the
Lie algebra complex and 
\begin{align}
	\delta:\Sym^1\g^\vee&\to \Sym^1\g^\vee\otimes \Lambda^{1}\g^\vee\\
	\nonumber	X^\vee&\mapsto \sum_{i=1}^k\theta(X_i)X^\vee\otimes h(X^\vee_i).
\end{align}
Here $\theta(X_i)X^\vee(Y):=X^\vee([Y,X_i])$. It is clear from this definition that 
\[
	W^{\ge n,\cdot}(\g):=\bigoplus_{p\ge n}W^{p,\cdot}(\g)
\]
is a subcomplex of $W^{*,\cdot}(\g)$. On the other hand,
\[ 
	W^{< n,\cdot}(\g):=\bigoplus_{p< n}W^{p,\cdot}(\g)
\]
 is a quotient of $W^{*,\cdot}(\g)$ and canonically isomorphic to $W^{*,\cdot}(\g)/W^{\ge n,\cdot}(\g)$.
For $n=1$ we have
\[
	W^{<1,\cdot}(\g)=\cC^\cdot(\g).
\]
The following lemma is classical
(see Cartan \cite{Ca} or \cite{Ra} Lemma 2.10. for a proof).
\begin{lemma}\label{weilalgcoh}
	One has:
	\begin{itemize}
	\item[a)] $H^0(W^{*,\cdot})=K$
	\item[b)] $H^n(W^{*,\cdot})=0$ for $n>0$
	\item[c)] $H^{2n}(W^{\ge n,\cdot})=(\Sym^n\g^\vee)^\g$
	\end{itemize}
\end{lemma}
Consider the exact sequence
\[
	0\to W^{\ge 1,\cdot}(\g)\to W^{*,\cdot}(\g)\to \cC^\cdot(\g) \to 0.
\]
This induces a connecting homomorphism
\[
	\partial:H^{2n-1}(\g)\to H^{2n}(W^{\ge 1,\cdot}(\g)),
\]
which is an isomorphism for $n>0$ by Lemma \ref{weilalgcoh}.
\begin{definition}\label{sfrgdefn}
	The {\em suspension} morphism for $\g$ is the composition 
	\[ 
		s_\g:H^{2n}(W^{\ge n,\cdot}(\g))\to H^{2n}(W^{\ge 1,\cdot}(\g))\xrightarrow{\partial^{-1}} H^{2n-1}(\g). 
	\]
\end{definition}
The suspension $s_\g$ has a different description:
\begin{lemma}\label{alternatesuspdefn}
	The following diagram commutes:
	\[
	\begin{xy}\xymatrix{
			H^{2n}(W^{\ge n,\cdot}(\g))\ar[rd]^{s_\g}\\
			H^{2n-1}(W^{< n,\cdot}(\g))\ar[r]\ar[u]^{\partial}& 
                               H^{2n-1}(\g).
	}\end{xy}
	\]
\end{lemma}
\begin{proof}
	This is just the compatibility of the two coboundary maps for the two short exact sequences.
\end{proof}

\subsection{The cosimplicial de Rham complex and the Weil algebra}\label{infdefn}
In this section we formulate the extension of Proposition \ref{normofAG} to the de Rham complex and
the Weil algebra. 


\begin{proposition}[\cite{Burgos} Proposition 8.10.]\label{normofAGinf}
	There is a natural bigraded isomorphism
	\[
		\cN \Omega^*(\BGinf^\alg)\isom W^{*,\cdot}(\frg).
	\]
\end{proposition}
The bigrading gives:
\begin{corollary}\label{normofAGinfcor}
	There are a natural isomorphisms 
	\begin{gather*}
		\cN \Omega^{\ge n}(\BGinf^\alg)\isom W^{\ge n,\cdot}(\frg),\\
		\cN \Omega^{< n}(\BGinf^\alg)\isom W^{< n,\cdot}(\frg).
	\end{gather*}
\end{corollary}
\begin{defn}\label{2infdef} We also denote by $\inf$ the map of complexes
induced by $\inf:\Ohalg{B.G}\to\Ohalg{\BGinf}$
\[ \cN\Omega^*(B.G^\alg)\to\cN\Omega^*(\BGinf^\alg)\isom W^{\ge n,\cdot}(\frg)\]
and the isomorphism of Corollary \ref{normofAGinfcor}.
\end{defn}
\begin{lemma}\label{diagr1} The map $\inf$ factors through $\cN\Omega^*(B.\cG^\dagger)$
	and the diagram
	\[
	\begin{CD}
		H^{2n}(\Omega^{\ge n}(\BcG^\dagger))@>{\inf}>>H^{2n}(W^{\ge n,\cdot}(\frg))\\
		@A\partial AA@AA\partial A\\
		H^{2n-1}(\Omega^{< n}(\BcG^\dagger))@>{\inf}>>H^{2n-1}(W^{< n,\cdot}(\frg))
	\end{CD}
	\]
	commutes, where the vertical maps are the boundary maps for the obvious 
	exact sequences.
\end{lemma}
\bew As $\inf:\Ohalg{B.G}\to\Ohalg{\BGinf}$ factors through $\Ohdagger{B.\cG}$ the first statement is clear.
	The morphism of short exact sequences 
	\[
	\begin{CD}
	\cN\Omega^{\ge n}(\BcG^\dagger)@>>>\cN\Omega^{*}(\BcG^\dagger)@>>>\cN\Omega^{< n}(\BcG^\dagger)\\
	@V{\inf}VV@V{\inf}VV@V{\inf}VV\\
	W^{\ge n,\cdot}(\frg)@>>>W^{*,\cdot}(\frg)@>>> W^{< n,\cdot}(\frg).
	\end{CD}
	\]
induces natural boundary maps.
\bewende

\subsection{Comparison of the suspension maps}
Now we can state the relation between the suspension $s_G$ for $B.G$ and its
infinitesimal version $s_\g$ on the Weil algebra. The proof is taken from Burgos \cite{Burgos}.

Consider $s_G$ as the composition  
\[
	s_G:H^{2n}(\Omega^{\ge n}(B.G^\alg))\to H^{2n-1}(\Omega^{\ge n}(G^\alg))\to H^{2n-1}_\DR(G^\alg)
\]
as in Section \ref{suspdefn}. 
\begin{prop}\label{suspensioncomp}
	There is a commutative diagram
	\[
	\begin{CD}
		H^{2n}(\Omega^{\ge n}(\BG^\alg))@>s_G>> H^{2n-1}_\DR(G^\alg)\\
		@V{\inf}VV@AA\rho A\\
		H^{2n}(W^{\ge n,\cdot}(\frg))@>s_\frg >>H^{2n-1}(\frg)
	\end{CD}
	\]	
	where $\rho$ is the isomorphism of Lemma \ref{rhodefn}.
\end{prop}
\begin{proof}
	Consider the map of complexes
	\[
		\cN \Omega^{\ge n}(B.G^\alg)\to \cN \Omega^{\ge n}(\BGinf^\alg)\isom W^{\ge n,\cdot}(\frg)
	\]
	defined in Definition \ref{2infdef}.
	According to \cite{Burgos} Theorem 8.12 $\inf$ induces a map 
	\[
		\omega_{E.G}^{-1}:H^{2n}(\Omega^{\ge n}(B.G^\alg)\to H^{2n}(W^{\ge n,\cdot}(\frg)),
	\]
	which is an algebraic description of the inverse of 
	the Chern-Weil homomorphism. Proposition 5.33 in \cite{Burgos} says that $s_G$ can be factored
	\[
		H^{2n}(\Omega^{\ge n}(B.G^\alg)\xrightarrow{\omega_{E.G}^{-1}}H^{2n}(W^{\ge n,\cdot}(\frg))\xrightarrow{s_\frg}
		H^{2n-1}(\frg).
	\]
	The proof in loc. cit. is only over $\bbC$, but works without any essential changes over an arbitrary field of
	characteristic zero.
	This gives the desired commutativity.
\end{proof}

\subsection{Proof of Theorem \ref{suspdiagram}}\label{proofofsuspensiondiag}

We
need one more  lemma.

\begin{lemma}\label{diagr2}
	The diagram 
	\[
	\begin{CD}
		H^{2n-1}(\Omega^{< n}(\BcG^\dagger))@>{\inf}>>H^{2n-1}(W^{< n,\cdot}(\frg))\\
		@VVV@VVV\\
		H^{2n-1}_\la(\cG,K)@>\Psi>>H^{2n-1}(\frg)
	\end{CD}
	\]
	commutes, where $\Psi$ is the map defined in Definition \ref{psidefn}.
\end{lemma}
\begin{proof}
	We have the commutative diagram
	\[
	\begin{CD}
		\cN\Omega^{< n}(\BcG^\dagger)@>{\inf}>>W^{< n,\cdot}(\frg)\\
		@VVV@VVV\\
		\cN\Omega^{< 1}(\BcG^\dagger)@>{\inf}>>W^{< 1,\cdot}(\frg)\\
		@VVV @VV=V\\
		\Ohdagger{\BcG}@>{\inf}>>\cC^\cdot(\frg).
	\end{CD}
	\]
By Lemma \ref{algebraicLazard} $\inf$ agrees with $\Psi$. Finally $\Psi$ factors naturally through $\Ohla{\BcG}$. 
\end{proof}

\begin{proof}[Proof of Theorem \ref{suspdiagram}]
We only have to combine the commutative diagrams that we have established:
The statements on $\inf$ were shown in Lemma \ref{diagr1}. The diagram for $s_G$ is Lemma \ref{suspensioncomp}. The small triangle was considered in Lemma \ref{alternatesuspdefn}. Finally the diagram for $\Psi$ is the above Lemma \ref{diagr2}.
\end{proof}


\section{The identification of $\Phi$ with the Lazard Isomorphism}

In this section we work with a smooth algebraic group $H/\Z_p$ and a
 $p$-saturated group of finite rank $\cH$ (see \cite{L} III Definition 2.1.3 and 2.1.6) with valuation $\omega$, which is an open subgroup of $H(\Z_p)$.
The group $\cH$ will always be considered as a $\bbQ_p$-analytic manifold.
Our main example is $H=\GL_N$, $\cH=1+pM_N(\Z_p)$ with $\omega$ the $\inf$-valuation.

\subsection{The Lazard Lie algebra}

Let $\Al\cH$ be the completed 
group ring of $\cH$ over $\bbZ_p$ (as defined in \cite{L} II 2.2.1.).
Note that by \cite{L} III. 3.3.2.1 $\cH$ is an analytic taylor manifold, in particular,
that $\cH$ can be identified with $\bbZ_p^r$ as an analytic manifold. Fix
$$
	\phi:\cH\to \bbZ_p^r \mbox{ with }\phi(e)=0
$$
such an identification.

If we are only interested in the structure of $\Al\cH$ as a $\bbZ_p$-module,
we can identify $\cH$ with $\bbZ_p^r$ via $\phi$ and obtain as in \cite{L} III 3.3.5
a topological basis $(z^\alpha)_{\alpha\in J}$, where $J:=\bbN^r$
(see \cite{L} III 2.3.8 and III 2.3.11.3). Here $z_i:=x_i-1$ for a ordered basis
$x_1,\dots, x_r$ of $\cH$ and $z^\alpha:=\prod_{i=1}^rz_i^{\alpha_i}$.
The valuation of $z^\alpha$ is
$$
	w(z^\alpha):=\sum_{i=1}^r\alpha_i\omega(x_i).
$$
This means that every element
$x\in \Al\cH$ can be written in the form
$$
	x=\sum_{\alpha \in J}\lambda_\alpha z^\alpha
$$
with $\lambda_\alpha\in\bbZ_p$
and the valuation is defined by
$$
	w(\sum_{\alpha \in J}\lambda_\alpha z^\alpha):=\inf_{\alpha\in J}\{v(\lambda_\alpha)+w(z^\alpha)\}
$$
(see \cite{L} I 2.1.17). The map $\cH\to \Al\cH$ is explicitly given by
$$
	\prod_{i=1}^rx_i^{\lambda_i}=\sum_{\alpha\in J}{\lambda\choose \alpha}z^\alpha.
$$
The saturation of $\Al\cH$ contains by definition the elements $\mu z^\alpha$, $\mu\in\bbQ_p$ with
$w(\mu z^\alpha)\ge 0$, i.e., with $v(\mu)+\sum_{i=1}^r\alpha_i\omega(x_i)\ge 0$.
As $\alpha_i\omega(x_i)>\frac{\alpha_i}{p-1}$ (see \cite{L} III 2.2.7.1) and
$v(\alpha_i!)\le \frac{\alpha_i}{p-1}$ we see that 
$$
	e_\alpha:=\frac{z^\alpha}{\alpha !}\in \Sat \Al\cH.
$$
The saturated group algebra $\Sat \Al\cH$ (see \cite{L} I 2.2.11) has also a structure 
of a valued, diagonal $\bbZ_p$-algebra, i.e., one has a valued $\bbZ_p$-algebra morphism
$$
	\Delta:\Sat\Al\cH\to \Sat\Al\cH\otimes_{\bbZ_p}\Sat\Al\cH
$$
of supplemented algebras. This is defined using the diagonal map
$\cH\to \cH\times \cH$ and \cite{L} I 3.2.8, II 2.2.8.  
\begin{definition}[\cite{L} IV 1.3.1]\label{LazardLie}
The \emph{Lazard Lie algebra} $\cL^*$ of $\cH$ is defined to
be
$$
	\cL^*:=\cL^* \Sat \Al\cH=\{x\in \Sat \Al\cH|\Delta(x)=x\otimes 1+1\otimes x\mbox{ and }
	w(x)>\frac{1}{p-1} \}.
$$
\end{definition}
We would like to make this more explicit. Let us first define special elements in $\Sat\Al\cH$:
\begin{definition}
	Let $\partial_i\in \Sat\Al\cH$ be the element
	$$
		\partial_i:=\sum_{\alpha_i> 0}\frac{(-1)^{\alpha_i-1}z_i^{\alpha_i}}{\alpha_i}=\log(z_i+1)=\log(x_i).
	$$
\end{definition}

\begin{lemma}
	The elements $\partial_i$ are primitive in $\Sat\Al\cH$, i.e., 
	$$
		\Delta(\partial_i)=\partial_i\otimes 1+1\otimes \partial_i.
	$$
	Moreover, they give a basis of $\cL^*$.
\end{lemma}
\begin{proof}
	This is just lemma IV 3.3.6 in \cite{L} as $z_i+1=x_i$.
\end{proof}
\begin{cor}
	There is a morphism
	$$
		\cU(\cL^*)\to \Sat \Al\cH.
	$$
\end{cor}
\begin{proof}
	Clear from the universal property of $\cU(\cL^*)$.
\end{proof}

\subsection{Distributions of locally analytic functions}
We show that the space of distributions of locally analytic functions
$\cD_\cont(\cH)$ is a subspace of $\Sat\Al\cH$.
It turns out that the elements $\partial_i\in \Sat\Al\cH$ can also be viewed as distributions of
locally analytic functions. 

Recall that a continuous function defined by the Mahler series
$$
	f(\lambda)=\sum_{\alpha\in J}c_\alpha {\lambda\choose \alpha}\mbox{ with $c_\alpha\in \bbZ_p$}
$$
is locally analytic if $\liminf_{|\alpha|\to \infty}\frac{v(c_\alpha)}{|\alpha|}>0$
(see \cite{L} III 1.3.9.2). 

Recall that
$\Ohla{\cH}\isom\Ohla{(\bbZ_p^r)}=\cup_{h>0}LA^h(\bbZ_p^r,\bbQ_p)$, where
$LA^h(\bbZ_p^r,\bbQ_p)$ are the locally analytic functions of order $h$
on $\bbZ_p^r$ with values in $\bbQ_p$ (see \cite{L} III 1.3.7). Each 
$LA^h(\bbZ_p^r,\bbQ_p)$ is a $p$-adic Banach space with norm $v_{LA^h}$
(see \cite{Co} 1.4.2 for the definition) and $\Ohla{\cH}$ gets the inverse
limit topology. We define
$$
	\cD_\cont(\cH):=\Hom_\cont(\Ohla{\cH},\bbQ_p).
$$
Amice shows in the case $r=1$, which extends immediately to
$r\ge 1$, the following proposition:
\begin{proposition}[Amice, cf. \cite{Co} Th\'eor\`eme 2.3]\label{Amice}
	The ring of distributions 
	$$	
	\cD_\cont(\cH)=\Hom_\cont(\Ohla{\cH},\bbQ_p)
	$$
	is identified via the \emph{Amice or Fourier transformation} 
	\begin{align*}
		\cA:\cD_\cont(\cH)&\to \bbQ_p[[T_1,\dots,T_r]]\\
		\mu&\mapsto \sum_\alpha\mu({\lambda\choose \alpha}) T^\alpha
	\end{align*}
	with the
	ring of power series $G(T)=\sum_\alpha\lambda_\alpha T^\alpha$, with $\lambda_\alpha\in\bbQ_p$,
	which converge for $v(T)>0$ (i.e., $ G(T)$ defines an analytic function on
	the open $r$-dimensional unit ball). 
\end{proposition}
\begin{cor}\label{Dcont}
	There is an injection
	\begin{align*}
		\cD_\cont(\cH)&\to \Sat\Al\cH\otimes_{\bbZ_p}\bbQ_p\\
		\mu&\mapsto \sum_\alpha\mu\left({\lambda\choose \alpha}\right)) z^\alpha.
	\end{align*}
	An element $D:=\sum_\alpha\rho_\alpha z^\alpha\in \Sat\Al\cH$ is an 
	element in $\cD_\cont(\cH)$ if and only if the power series
	$G_D(T):=\sum_\alpha{\rho_\alpha}T^\alpha$ converges for $v(T)>0$.
	Explicitly, for $f(\lambda)=\sum_{\alpha\in J}c_\alpha {\lambda\choose \alpha} \in\Ohla{\cH}$
	one has
	$$
		Df=\sum_\alpha c_\alpha\rho_\alpha.
	$$
\end{cor}

\begin{proof}
	Using Proposition \ref{Amice} we see that $\sum_\alpha\mu({\lambda\choose \alpha}) z^\alpha$
	converges in $\Sat\Al\cH$ as $w(z_i)=\omega(x_i)>\frac{1}{p-1}$ and that
	the map is injective. The rest of the corollary is just a restatement of
	Proposition \ref{Amice}.
\end{proof}
As $\partial_i=\sum_{\alpha_i> 0}\frac{(-1)^{\alpha_i-1}}{\alpha_i}z_i^{\alpha_i}$
we get:
\begin{cor}\label{derivativesofla}
	The elements $\partial_i$ are contained in $\cD_\cont(\cH)$ and the
	inclusion $\cL^*\subset \cD_\cont(\cH)$ defines a ring homomorphism
	$$
		\cU(\cL^*)\to \cD_\cont(\cH).
	$$
\end{cor}

\subsection{Identification of the algebraic Lie algebra $\frh$ with the Lazard Lie algebra $\cL^*$}
\label{ident}
Recall that the algebraic $\bbQ_p$-Lie algebra of a smooth linear algebraic group scheme $H/\bbZ_p$ is defined in Definition \ref{lie-algebra} as the derivations
$$
	\frh:=\Der_{\bbZ_p}(\Ohalg{H}_e,\bbQ_p).
$$
Inside $H(\bbZ_p)$ we have the open, $p$-saturated subgroup $\cH$ considered as $\bbQ_p$-analytic manifold.
Using the inclusion $\Ohalg{H}_e\subset \Ohla{\cH}_e\subset\wh{\Ohalg{H}_e}$ (completion w.r.t. the maximal ideal),
we can compare this with the Lazard Lie algebra $\cL^*$. 

\begin{prop}\label{liealgident}
	The inclusion $\Ohalg{H}_e\subset \Ohla{\cH}_e$ induces an isomorphism
	$$
		\frh\isom \cL^*\otimes_{\bbZ_p}\bbQ_p.
	$$
	Moreover, let
	$f(\lambda)=\sum_{\alpha\in J}c_\alpha{\lambda\choose \alpha}$ be a locally analytic function and $\partial_i\in \cL^*$, then
	$$
		\partial_if=\frac{\partial f}{\partial \lambda_i}(e).
	$$
\end{prop}
\begin{proof}
	We compute
	$$
		\frac{\partial f}{\partial \lambda_i}=\sum_{\alpha\in J}c_\alpha\frac{\partial }{\partial \lambda_i}{\lambda\choose \alpha}
	$$
	and 
	$$
		\frac{\partial }{\partial \lambda_i}{\lambda\choose \alpha}|_{\lambda=0}=\left\{\begin{array}{ll}
		\frac{(-1)^{\alpha_i-1}}{\alpha_i}\mbox{ if }\alpha=(0,\cdots, \alpha_i,\cdots, 0)\\
		0\mbox{ else}
		\end{array}\right..
	$$
	Thus, 
	$$
		\frac{\partial f}{\partial \lambda_i}|_{\lambda=0}=\sum_{\alpha_i>0}
		\frac{(-1)^{\alpha_i-1}c_{\alpha_i}}{\alpha_i}=\partial_if,
	$$
	which gives the desired result. To prove that $\frh\isom\cL^*\otimes_{\bbZ_p}\bbQ_p$, first observe
	that the maximal ideal of $\Ohla{\cH}_e$ is $\frm_e\Ohalg{G}_e$ and that the 
	$\partial_i$ by the above calculation form a basis of this ideal. As $\cL^*$ is generated by the
	$\partial_i$ we get that $\cL^*\otimes_RK\isom \frm_e\Ohalg{G}_e/\frm_e^2\Ohalg{G}_e$.
	As $\frh=\frm_e/\frm_e^2$ this proves that $\frh\isom\cL^*\otimes_{\bbZ_p}\bbQ_p$.
\end{proof}

\subsection{Standard complexes for group and Lie algebra cohomology}

For any augmented $\bbZ_p$-algebra $A$ with augmentation $\epsilon:A\to\bbZ_p$
we consider two \emph{standard complexes} $(T.A,d)$ and $(\wt T.A,\wt d)$ with 
$$
	T_nA=\wt T_nA=A^{\otimes n+1},
$$
where for simplicity
\begin{align*}
	A^{\otimes n+1}=\begin{cases}
				A\otimes_{\bbZ_p}\ldots\otimes_{\bbZ_p}A&\mbox{ $n+1$-times; if $A$ has no topology}\\
				A\wh\otimes_{\bbZ_p}\ldots\wh\otimes_{\bbZ_p}A&\mbox{ $n+1$-times; if $A$ complete}.
	                \end{cases}
\end{align*}
and differentials 
\begin{align*}\label{Lazarddifferential}
	d(a_0\otimes\ldots\otimes a_n):=&
		\sum_{i=0}^{n}(-1)^i\epsilon(a_i)a_0\otimes\ldots\otimes \wh a_i\otimes\ldots \otimes a_n\\
	\wt d(a_0\otimes\ldots\otimes a_n):=&\sum_{i=0}^{n-1}(-1)^ia_0\otimes\ldots\otimes a_ia_{i+1}\otimes\ldots \otimes a_n\\ &+(-1)^na_0\otimes\ldots\otimes a_{n-1}\epsilon(a_n)	
\end{align*}
To motivate these constructions, we consider two realizations $E.\cH$ and $\wt E.\cH$ of the
universal bundle over the classifying space $B.\cH$.
Let $E_n\cH=\wt E.\cH=\cH^{n+1}$ with $\delta_i(h_0,\dots,h_n)=(h_0,\dots,\hat h_i,\dots, h_n)$
and 
$$
	\wt\delta_i(h_0,\dots,h_n)=\begin{cases}
					(h_1,\dots,h_n)&\mbox{ if $i=0$}\\
					(h_0,\dots,h_{i+1}h_{i+2},\dots,h_n)&\mbox{ if $i>0$}
	                           \end{cases}
$$
The classifying space $B. \cH$ is then the quotient of these spaces by the $\cH$-action:
diagonally on $E.\cH$ and on the last factor on $\wt E.\cH$.

Associated to these contractible simplicial spaces we have simplicial
complexes $\cD_\cont(E.\cH), \Al(E.\cH), \Sat\Al(E.\cH)$ and 
similarly for $\wt E.\cH$ with differential the alternating sums of the
$\delta_i$'s respectively the $\wt \delta_i$'s.
\begin{lemma}
	One has
	$$
		\cD_\cont(E.\cH)\isom T.\cD_\cont(\cH)
	$$
	and
	$$
		\cD_\cont(\wt E.\cH)\isom\wt T.\cD_\cont(\cH)
	$$
	and similar results for $\Al\cH$ and $\Sat\Al\cH$.
\end{lemma}
\begin{proof}
	Clear from the definition.
\end{proof}
>From the isomorphism of simplicial spaces $E.\cH\isom \wt E.\cH$ one sees that there
is an isomorphism of complexes $\cD_\cont(E.\cH)\isom \cD_\cont(\wt E.\cH)$.
In a similar way one sees:
\begin{lemma}\label{cUisom}
	There is an isomorphism of complexes
	$$
		T.\cU\cL^*\isom \wt T.\cU\cL^*
	$$
	and both complexes are projective resolutions of $\bbQ_p$ as trivial $\cU\cL^*$-module.
\end{lemma}

\subsection{Review of the Lazard isomorphism}
Recall that $\cH$ is a $p$-saturated group of finite rank with valuation $\omega$
and $\cL^*$ is its Lazard Lie algebra.
The main theorem of chapter V in \cite{L} can be formulated as follows:
\begin{theorem}[\cite{L} V 2.4.9]\label{lazardisom}
	There is an isomorphism
	$$
		H^i_\cont(\cH,\bbQ_p)\isom H^i(\cL^*,\bbQ_p).
	$$
\end{theorem}
Let us review how Lazard constructs this isomorphism. First Lazard \cite{L} V 1.2.9 shows that
$\Hom_{\Al\cH}(\wt T.\Al\cH,\bbQ_p)$ computes the continuous group cohomology. 
Then the isomorphism is obtained from the following three quasi-isomorphisms:
\begin{equation}\label{qis1}
	\Hom_{\Al\cH}(\wt T.\Sat\Al\cH,\bbQ_p)\xrightarrow{qis}\Hom_{\Al\cH}(\wt T.\Al\cH,\bbQ_p)
\end{equation}
induced from $\Al\cH\subset \Sat\Al\cH$, 
\begin{equation}\label{qis2}
	\Hom_{\Sat\Al\cH}(\wt T.\Sat\Al\cH,\bbQ_p)
	\xrightarrow{qis}\Hom_{\cU\cL^*}(\wt T.\cU\cL^*,\bbQ_p),
\end{equation}
induced from $\cU\cL^*\to \Sat\Al\cH$ and finally
\begin{equation}\label{qis3}
	\Hom_{\cU\cL^*}(\wt T.\cU\cL^*,\bbQ_p)\xrightarrow{qis}\Hom_{\cU\cL^*}(\cU\cL^*\otimes{\bigwedge}^\cdot\cL^*,\bbQ_p)
\end{equation}
induced from the anti-symmetrisation map $as_n$ 
$$
	as_n:{\bigwedge}^n\cL^*\to \cU(\cL^*)^{\otimes n}
$$
given by
$$
	as_n(X_1\wedge\ldots\wedge X_n)=\sum_{\sigma\in \cS_n}\sgn(\sigma)X_{\sigma^{-1}(1)}\otimes
	\ldots\otimes X_{\sigma^{-1}(n)}.
$$
The fact that the latter is a quasi-isomorphism follows from \cite{CE} XIII Theorem 7.1.

\subsection{Explicit description of the Lazard isomorphism}
We describe the Lazard isomorphism as a kind of Taylor series expansion.

In the last section the Lazard isomorphism was shown to be induced from the map of 
complexes 
$$
	\cU\cL^*\otimes{\bigwedge}^\cdot\cL^*\to \wt T.\cU\cL^*\to \wt T.\Sat\Al\cH.
$$
In Corollary \ref{Dcont} we saw that the map $\cU\cL^*\to \Sat\Al\cH\otimes_{bbZ_p}\bbQ_p$ factors 
through $\cD_\cont(\cH)$. We get a commutative diagram
$$
\begin{xy}\xymatrix{\Hom_{\Sat\Al\cH}(\wt T.\Sat\Al\cH,\bbQ_p)\ar[r]\ar[d]&
	\Hom_{\cU\cL^*}(\cU\cL^*\otimes{\bigwedge}^\cdot\cL^*,\bbQ_p)\\
	\Hom_{\cD_\cont(\cH)}(\wt T.\cD_\cont(\cH),\bbQ_p)\ar[ur]&}
	\end{xy}
$$
Using the identification
$$
	\Hom_{\cD_\cont(\cH)}(\wt T_n\cD_\cont(\cH),\bbQ_p)
	\isom \Hom_{\cont}(\cD_\cont(\cH)^{\otimes n},\bbQ_p)\isom \Ohla{\cH}^{\otimes n}
$$
and 
$$
	\Hom_{\cU\cL^*}(\cU\cL^*\otimes{\bigwedge}^\cdot\cL^*,\bbQ_p)\isom
	\Hom_{\bbZ_p}({\bigwedge}^\cdot\cL^*,\bbQ_p)
$$
the diagram gives:
$$
	\begin{xy}\xymatrix{\Hom_{\Sat\Al\cH}(\wt T.\Sat\Al\cH,\bbQ_p)\ar[r]\ar[d]&
	\Hom_{\bbZ_p}({\bigwedge}^\cdot\cL^*,\bbQ_p)\\
	\Ohla{\cH}^{\otimes n}\ar[ur]&}
	\end{xy}
$$
We want to make the map $\Ohla{\cH}^{\otimes n}\to \Hom_{\bbZ_p}({\bigwedge}^n\cL^*,\bbQ_p)$
more explicit.

Recall that we have defined in Section \ref{lie}
	for each $f\in \Ohla{\cH}$ a  linear form 
	$$
	df(e)\in \Hom_{\bbZ_p}(\cL^*,\bbQ_p)\isom \Hom_{\bbQ_p}(\frh,\bbQ_p)=\frh^\vee,
	$$ 
where the isomorphism comes from Proposition \ref{liealgident}. One has
	$$
 		df(e)(\partial_j)=\frac{\partial f}{\partial\lambda_j}(e).
 	$$
	This is the differential of $f$ evaluated at $e\in\cH$.
\begin{prop}\label{taylorent}
	The map
	$$
		\Ohla{\cH}^{\otimes n}\to \Hom_{\bbZ_p}({\bigwedge}^n\cL^*,\bbQ_p)\isom {\bigwedge}^n\frh^\vee
	$$
	is given by
	$$
	f_1\otimes\ldots\otimes f_n\mapsto df_1(e)\wedge\ldots \wedge df_n(e).
	$$
	Thus the Lazard isomorphism agrees with $\Phi$ as defined in Definition \ref{Phidefn}.
\end{prop}
\begin{proof}
	By definition, $f_1\otimes\ldots\otimes f_n$ maps to the linear form, which maps
	$\partial_1\wedge\ldots\wedge \partial_n\in {\bigwedge}^n\frh$ to
	$$
		\sum_{\sigma\in\cS_n}\sgn(\sigma)\partial_{\sigma^{-1}(1)}f_1
		\ldots \partial_{\sigma^{-1}(n)}f_n.
	$$
	But this is exactly the linear form $df_1(e)\wedge\ldots \wedge df_n(e)$  with
	$df_i(e)$ as defined above. 
\end{proof}
\begin{remark}
	It is not hard to see that the map 
	$$
		\Ohla{\cH}^{\otimes n}\to \Hom_{\bbZ_p}(\cU(\cL^*)^{\otimes n},\bbQ_p)
	$$
	is given by the whole Taylor series and not just the first coefficient. For our
	purposes the above result suffices.
\end{remark}
\subsection{Comparison of $\Phi$ with $\Psi$}
We return to the situation in Section~\ref{ident}, where we had a smooth algebraic group scheme $H/\bbZ_p$
and $\cH$ was a $p$-saturated open subgroup of $H(\bbZ_p)$.
Using Proposition \ref{liealgident} we identify $\cL^*\otimes\bbQ_p\isom \frh$.

In this section we relate the map of complexes
\begin{align*}
	\Phi:\Ohla{B_n\cH}&\to \cC^n(\frh)\\
	f_1\otimes\ldots\otimes f_n&\mapsto df_1(e)\wedge\ldots\wedge df_n(e)
\end{align*}
to the map of complexes
\begin{align*}
	\Psi:\Ohla{B_n\cH}\subset \Ohla{E_n\cH}&\to \cC^n(\frh)\\
	f_0\otimes\ldots\otimes f_n&\mapsto f_0(e)df_1(e)\wedge\ldots\wedge df_n(e)
\end{align*}
defined in Definition \ref{psidefn} (for $\cH=1+pM_N(\Z_p)$). We have the following theorem:
\begin{theorem}\label{lazgleichpsi}
	The map $\Phi$ and the map $\Psi$ are homotopic maps of complexes. In particular, they
	induce the same map on cohomology
	$$
		\Phi=\Psi:H^i_\la(\cH,\bbQ_p)\isom H^i(\frh,\bbQ_p).
	$$
\end{theorem}
\begin{proof}
	Write $\cC^n(\frh)=\Hom_{\cU\frh}(\cU\frh\otimes{\bigwedge}^n\frh,\bbQ_p)$, then the map 
	$\Ohla{E_n\cH}\to \cC^n(\frh)$, which defines $\Psi$ is induced by a  map of complexes
	$$
		\cU\frh\otimes{\bigwedge}^\cdot\frh\to T_\cdot\cU\frh\to T_\cdot\cD_\cont(\cH).
	$$
	On the other hand, $\Phi:\Ohla{B_n\cH}\to \cC^n(\frh)$ is induced in degree $n$ by a map
	$$
		{\bigwedge}^n\frh\to \cU\frh^{\otimes n }\to \cD_\cont(\cH)^{\otimes n}.
	$$
	We extend the first map $\cU\frh$-linearly to a map of complexes
	$$
		\cU\frh\otimes{\bigwedge}^\cdot\frh\to \wt T_\cdot\cU\frh.
	$$
	Using the isomorphism $ \wt T_\cdot\cU\frh\isom T_\cdot\cU\frh$ from Lemma $\ref{cUisom}$, 
	we get two maps from the projective resolution 
	$\cU\frh\otimes{\bigwedge}^\cdot\frh$ of $\bbQ_p$ to the resolution $T_\cdot\cU\frh$, which must be homotopic 
	by general facts for projective resolutions. Composing this with the commutative diagram
	$$
		\begin{CD}
			T.\cU\frh@>>> T.\cD_\cont(\cH)\\
			@VV\isom V@VV\isom V\\
			\wt T.\cU\frh@>>> \wt T.\cD_\cont(\cH)\\
		\end{CD}
	$$
	gives still homotopic maps, so that $\Phi$ and $\Psi$ coincide on cohomology.
\end{proof}
\section{Proof of the Main Theorem}\label{proofsection}
In this section we are going to put together the proofs of the results announced in Chapter \ref{setup}.

Let throughout $G=\GL_N$ as algebraic group over $R$, $\cG=1+\pi M_N(R)$ as
$K$-Lie group and $\cG^\dagger$ the underlying dagger-space.

\begin{proof}[Proof of Theorem \ref{lazardiso}]
Let $\cG_n=1+\pi^nM_N(R)$. They are open and closed normal subgroups of $\GL_N(R)$ of finite index. 
They are also
a neighbourhood basis of $e\in\GL_N(R)$.
Hence $\Ohla{G}_e$, the ring of germs of analytic functions, is given by $\lim\limits_{n\to \infty} \Oh(\cG_n^\la)$.
 For $n'>n$ the natural restriction maps
\[ H^i_\la(\cG_n,K)\to H^i_\la(\cG_{n'},K)\]
are isomorphisms as the we working with rational coefficients. Passing to the limit
\[ H^i_\la(\cG_n,K)\to H^i(\Ohla{\BG}_e)\]
is also an isomorphism. We define
\[ \Phi:H^i(\Ohla{\BG}_e)\to H^i(\g,K)\]
by 
\[ f_1\otimes\ldots\otimes f_n\mapsto df_1(e)\wedge\ldots \wedge df_n(e).\]
Note that this is the same formula as in Definition \ref{Phidefn}.

Now let $R=\Z_p$. Then $\Phi$ is an isomorphism by Theorem \ref{lazardisom} together with Proposition \ref{taylorent}. Moreover, 
\begin{gather*}
\Lie( \GL_N(R))\isom \Lie(\GL_N(\Z_p))\tensor K\\
H^i(\Lie(\GL_N(R)),K)\isom H^i(\Lie(\GL_N(\Z_p)),\Q_p)\tensor K\\
\Oh(G_R^\la)_e\isom \Oh(G_{\Z_p}^\la)_e\tensor K
\end{gather*}
and $\Phi$ is compatible with extension of scalars. Hence  $\Phi$ is also an isomorphism for general $K$.
\end{proof}

\begin{proof}[Proof of Theorem \ref{mainthm}]
Let $1<n\leq N$. Putting together the commutative diagrams of Proposition \ref{etale}, Proposition \ref{diagrams} and Theorem \ref{suspdiagram}, we have established
the following big commutative diagram:

{
\footnotesize
\[\begin{xy}\hbox to -1.5cm{}
\xymatrix{
  &&H^{2n}(\Omega^{\geq n} \BG^\alg) \ar[d]\ar[d]\ar[rr]^{s_G}&&H^{2n-1}_\DR(G^\alg)\ar[ddd]^\rho\\
H^{2n}_\et(\BG,n)\ar[d]&&H^{2n}(\Omega^{\geq n}\BcG^\dagger)\ar[r]^\inf &H^{2n}(W^{\geq n,*}(\g))\ar[ddr]^{s_\g}\\
H^{2n-1}(G(R),H^1_\et(K,n))&H^{2n}_\syn(\BG,n)\ar[ul]\ar[d]\ar[r]^\eta &H^{2n-1}(\Omega^{<n}\BcG^\dagger)
       \ar[r]\ar[u]^\partial\ar[d]&H^{2n-1}(W^{<n,*}(\g))\ar[dr]\ar[u]\\
&H^{2n-1}(\cG,K)\ar[ul]^{\exp_{BK}} &H^{2n-1}_\la(\cG,K)\ar[l]\ar^{\Psi=\Phi}[rr]&&H^{2n-1}(\g,K)
}\end{xy}
\]
}

By Theorem \ref{lazgleichpsi} $\Psi=\Phi$ is the Lazard isomorphism. Note also that $H^{2n-1}(\cG,K)\isom H^{2n-1}(G(R),K)$ because $\cG$ is of finite index in $G(R)$.

We start with the primitive element $p_n\in H^{2n-1}(\g,K)$ and want
to identify its image in $H^{2n-1}(\cG,K)$. The Chern class $c_n^\alg\in H^{2n}(\Omega^{\geq n}\BG^\alg)$ is mapped to $p_n$ under the suspension map $s_G$. Using Proposition \ref{suspensioncomp},
\cite{Greub} VI.6.19. (which computes the image of the Chern class under $s_\g$) and \cite{Burgos} Lemma 8.11, one verifies the normalization of $p_n$ used
in Definition \ref{primitive}. 
By Proposition \ref{diagrams} the image of $c_n^\alg$ in $H^{2n}(\Omega^{\geq n}\BcG^\dagger)$ agrees with the image of the syntomic Chern class
$c_n^\syn\in H^{2n}_\syn(\BG,n)$ under $\partial\verk\eta$. By the commutativity of the
diagram this implies that the images of $p_n$ and $c_n^\syn$ in $H^{2n-1}(\cG,K)$ agree.
By Proposition \ref{etale} the universal syntomic Chern class is mapped to the universal \'etale Chern class. Together this proves the theorem. 
The case $n=1$ uses
basically the same argument.
\end{proof}


\begin{thebibliography}{99999}
\bibitem[Ba]{Ba}K. Bannai, Syntomic cohomology as a $p$-adic absolute Hodge cohomology, Math. Z.  242  (2002),  no. 3, 443--480.
\bibitem[B]{Bei}A. Beilinson, Higher regulators and values of $L$-functions, J. Soviet Math. {\bf 30} (1985), 2036-2070.
\bibitem[Be]{Be}A. Besser, Syntomic regulators and $p$-adic integration I: rigid syntomic regulators, Israel J. of Math. 120 (2000), 291--334.
\bibitem[Ber]{Ber}P. Berthelot, Cohomologie rigide et cohomologie rigide \`a supports propres, Pr\'epublication 96-03 Universit\'e de Rennes, 1996.
\bibitem[Bor1]{Bor}A.\ Borel,
Stable real cohomology of arithmetic groups,
Ann. Sci. \'Ecole Norm. Sup. (4) 7 (1974), 235--272.
\bibitem[Bor2]{Bo}A.\ Borel, Cohomologie de ${\rm SL}\sb{n}$ et valeurs de fonctions zeta aux points entiers, Ann. Scuola Norm. Sup. Pisa Cl. Sci. (4)  4  (1977), no. 4, 613--636.
\bibitem[BGR]{BGR}
S. Bosch, U. G\"untzer, R Remmert, 
Non-Archimedean analysis.
A systematic approach to rigid analytic geometry. Grundlehren der Mathematischen Wissenschaften [Fundamental Principles of Mathematical Sciences], 261.
Springer-Verlag, Berlin, 1984. 
\bibitem[BK]{BK}S. Bloch, K. Kato, $L$-functions and Tamagawa numbers of motives.  The Grothendieck Festschrift, Vol. I,  333--400, Progr. Math., 86, Birkh\"auser Boston, Boston, MA, 1990. 
\bibitem[Bu]{Burgos}J. I. Burgos Gil, The Regulators of Beilinson and Borel, CRM Monograph Series, 15. American Mathematical Society, Providence, RI, 2002.
\bibitem[Ca]{Ca}H. Cartan: Notions d'alg\`ebre diff\'erentielle, \OE uvres, vol. III
(R. Remmert and J.-P. Serre, eds.), Springer Verlag, 1979, 1268-1282.
\bibitem[CE]{CE}H. Cartan, S. Eilenberg: Homological algebra, Princeton University Press, 1956
\bibitem[Co]{Co}P. Colmez: Fonctions d'une variable p-adique, Preprint 2005 available at
www.math.jussieu.fr/~colmez/publications.html
\bibitem[FM]{FM}J.-M. Fontaine, W. Messing, 
$p$-adic periods and $p$-adic �ale cohomology.
Current trends in arithmetical algebraic geometry (Arcata, Calif., 1985), 179--207,
Contemp. Math., 67,
Amer. Math. Soc., Providence, RI, 1987. 
\bibitem[GHV]{Greub}W. Greub, S. Halperin, R. Vanstone, Connections, curvature and cohomology III,
Academic Press, 1976.
\bibitem[G]{G}M. Gros, R\'egulateurs syntomiques et valeurs de fonctions $L$ $p$-adiques II, Invent. math. 115, 61-79 (1994).
\bibitem[GK1]{GK}E. Gro\ss{}e-Kl\"onne, de Rham-Kohomologie in der rigiden Analysis, Dissertation M\"unster 1999.
\bibitem[GK2]{GK2}E. Gro\ss{}e-Kl\"onne, Rigid analytic spaces with overconvergent structure sheaf, J. reine angew. Math. 519 (2000), 73--95.
\bibitem[GK3]{GK3}E. Gro\ss{}e-Kl\"onne, Finiteness of de Rham cohomology in rigid analysis, Duke Math. J. 113, No. 1 (2002), 57--91.
\bibitem[Ha]{Hamida}N. Hamida: Les r\'egulateurs en $K$-th\'eorie alg\'ebrique, These Jussieu, 2002.
\bibitem[Ho]{Ho}G. Hochschild, Cohomology of algebraic linear groups,  Illinois J. Math.  5  1961 492--519.
\bibitem[Kar]{Karoubi}M. Karoubi, Sur la $K$-th\'eorie Multiplicative, in: J. Cuntz, M. Khalkali (eds.),
Cyclic Cohomology and Noncommutative Geometry, Fields Institute Communications 1997. 
\bibitem[Kat]{K}K. Kato, On $p$-adic Vanishing Cycles (Application  of Ideas of Fontaine-Messing),  Proc. Algebraic geometry, Sendai, 1985,  207--251, Adv. Stud. Pure Math., 10, North-Holland, Amsterdam, 1987.
\bibitem[L]{L}M. Lazard, Groupes analytiques $p$-adiques, Publ. IHES No. 2
6, 1965.
\bibitem[Lod]{Lod}J.-L. Loday,  Cyclic homology.  Grundlehren der Mathematischen Wissenschaften, 301. Springer-Verlag, Berlin, 1992.
\bibitem[Ni1]{Niz}W. Niziol, Cohomology of crystalline smooth sheaves, Compositio Math.  129  (2001),  no. 2, 123--147.
\bibitem[Ni2]{Ni2}W. Niziol, On the image of $p$-adic regulators, Invent. math. 127, 375--400 (1997).
\bibitem[Ra]{Ra}M. Rapoport: Comparison of the regulators of Beilinson and Borel, in: Be\u\i linson's conjectures on special values of $L$-functions,  169--192,
Perspect. Math., 4, Academic Press, Boston, MA, 1988. 
\bibitem[Se]{Serre}J.-P. Serre: Lie algebras and Lie groups, W. A. Benjamin, New York-Amsterdam 1965.
\bibitem[So1]{Soule}C. Soul\'e, $K$-th\'eorie des anneaux d'entiers de corps de nombres et cohomologie \'etale, Invent. Math. {\bf 55},  (1979), 251--295.
\bibitem[So2]{So2}C. Soul\'e, On higher $p$-adic regulators.  Algebraic $K$-theory, Evanston 1980 (Proc. Conf., Northwestern Univ., Evanston, Ill., 1980),  pp. 372--401, Lecture Notes in Math., 854, Springer, Berlin-New York, 1981. 
\end{thebibliography}
\end{document}